%% file: autosam.tex
\documentclass[twocolumn]{autart}    

\input{Files/preamble}

\newcommand{\lhk}[1]{\textcolor{black}{#1}} 
\newcommand{\lhl}[1]{\textcolor{black}{#1}} 
\newcommand{\ils}[1]{\textcolor{black}{#1}} 
\newcommand{\ilc}[1]{\textcolor{black}{#1}} 
\newcommand{\lhn}[1]{\textcolor{black}{#1}} 
\newcommand{\lho}[1]{\textcolor{black}{#1}} 
\newcommand{\ilo}[1]{\textcolor{black}{#1}} 
\begin{document}

\begin{frontmatter}

\title{An Optimal Control Interpretation of Augmented Distributed Optimization Algorithms\thanksref{footnoteinfo}} 

\thanks[footnoteinfo]{\lhk{
\ils{
This manuscript is a significantly extended version of
conference paper \cite{hallinan_InverseOptimal_24},
which includes the proofs of the results, additional discussion, and generalizations to wider classes of distributed optimization algorithms involving coupling inequality constraints and feed-forward terms.}}}

\author[Cambridge]{Liam Hallinan}\ead{lh706@cam.ac.uk},    
\author[Cambridge]{Ioannis Lestas}\ead{icl20@cam.ac.uk}               

\address[Cambridge]{Department of Engineering, University of Cambridge, Trumpington Street, Cambridge, CB2 1PZ, United Kingdom}  

\begin{abstract}                          
\li{Distributed optimization}
algorithms are used in a wide variety of problems
\lhe{involving complex network systems}
where the goal is for a set of agents in the network to solve a network-wide
\li{optimization} problem via \li{distributed} 
update rules. \lhe{In \il{many applications}, such as communication networks and power systems, transient performance of the algorithms is just as critical as convergence, as the algorithms link to physical processes which must behave well. \li{Primal-dual algorithms have a long history in solving distributed optimization problems, with augmented Lagrangian methods leading to important classes of widely used algorithms,  which have been observed 
 in simulations to improve transient performance.
Here we show that 
such algorithms can be seen as being the optimal solution to an appropriately formulated optimal control problem, i.e., a cost functional associated with the transient \lhi{behavior} of the algorithm is minimized, penalizing deviations from optimality during algorithm transients. This is shown for broad classes of algorithm dynamics, including the more involved setting where inequality constraints are present.
The results presented improve our understanding of \lic{the} performance of distributed optimization algorithms and can be used as a basis for improved formulations.}}
%
%
\end{abstract}

\end{frontmatter}

\section{Introduction} \label{sec_introduction}
The field of distributed optimization has gained increasing prominence in recent years in response to the increasing scope and complexity of many large-scale decentralized systems, such as the internet and power systems.
A central pillar in this field is the primal-dual framework for convex optimization, originally motivated by mid-twentieth-century economic theory, most notably within the
\lhk{resource-allocation problems studied by} Arrow, Debreu, Hurwicz, and Uzawa \cite{arrow_ExistenceEquilibrium_54,arrow_StabilityCompetitive_58, arrow_StudiesLinear_58}.
In these methods, distributed dual variables are adjusted to ensure constraint satisfaction, corresponding to a pricing mechanism, while local primal variables are updated via gradient-based rules and dual feedback \cite{boyd_ConvexOptimization_04, bertsekas_NonlinearProgramming_16}. The resulting dynamics can equivalently be viewed as seeking the saddle point of the Lagrangian function describing the optimization problem, where the dual variables correspond to Lagrange multipliers, combined with a non-smooth projection operator to ensure satisfaction of inequality constraints \cite{fletcher_PracticalMethods_00, nedic_SubgradientMethods_09}.

In network systems where decentralized agents collectively solve a global objective, these algorithms naturally decouple into distributed dynamics, making them particularly well suited.
In communication networks, for example, dual variables capture link congestion, enabling local flow-control decisions without a central coordinator \cite{kelly_RateControl_98, srikant_MathematicsInternet_04, low_AnalyticalMethods_22}. Similarly, power \lhj{systems} can use distributed methods to dispatch generators or loads to meet economic and operational constraints \cite{molzahn_SurveyDistributed_17,zhao_DesignStability_14,kasis_PrimaryFrequency_17,stegink_UnifyingEnergyBased_17}. Additional applications include machine learning \cite{boyd_DistributedOptimization_11,dean_LargeScale_12}, distributed parameter estimation \cite{xiao_DistributedAverage_07, kar_DistributedParameter_12}, robot formation control \cite{ren_DistributedConsensus_08, freeman_DistributedEstimation_06}, transportation \cite{alonso-mora_OndemandHighcapacity_17}, and opinion dynamics \cite{golub_NaiveLearning_10}.

Primal–dual algorithms often perform poorly when constraints are ill-conditioned, motivating the use of augmented Lagrangian techniques that add penalty or regularization terms to the system Lagrangian. \lhh{Such techniques include the Hestenes-Powell (or multiplier) method \cite{hestenes_MultiplierGradient_69, powell_MethodNonlinear_69, rockafellar_AugmentedLagrange_74, bertsekas_ConstrainedOptimization_82}, the alternating direction method of multipliers (ADMM) \cite{boyd_DistributedOptimization_11}, and proximal-point algorithms \cite{rockafellar_AugmentedLagrangians_76}, which have been generalized by considering more generic additions to the Lagrangian \cite{gill_PrimaldualAugmented_12, yamashita_PassivitybasedGeneralization_20}. Applications to network problems have been reviewed in \cite{jakovetic_PrimalDual_20}, and include the proportional-integral (PI) distributed optimization algorithm \cite{wang_ControlPerspective_11, droge_ContinuoustimeProportionalintegral_14}}. These modified algorithms seek the saddle point of the augmented Lagrangian and have been shown via empirical studies to exhibit superior performance to their standard counterparts \cite{birgin_NumericalComparison_05, holding_StabilityInstability_21a, holding_StabilityInstability_21, wang_ControlPerspective_11, feijer_StabilityPrimal_10}.

In the optimization literature, performance is typically benchmarked by convergence rate \cite{bertsekas_ConstrainedOptimization_82, ruszczynski_ConvergenceAugmented_95, nedic_SubgradientMethods_09, attouch_ConvergenceRate_20}. However, such metrics can only provide loose upper bounds \lhh{on settling time
(e.g., linearly in the running time \cite{zhang_DistributedOptimization_18}, which can be related to the spectrum of the underlying network graph \cite{jakovetic_LinearConvergence_15})
that are not easily tailored to specific problems},
and often overlook transient \lhj{behaviors} (e.g., overshoot or oscillations). These transient properties are particularly critical in real engineering systems as they correspond to physical \lhj{processes} -- such as power generation or data transmission.
\lhj{Analyzing} the transient \lhj{behavior} for such systems is challenging in general due to the nonlinear, non-smooth and distributed nature of the underlying dynamics, meaning the algorithm parameters are typically designed using heuristic techniques.

\lhh{To address these issues}, we propose an optimal control approach.
\lhk{In this framework, augmented algorithms can be interpreted as dynamical control systems (similar to \cite{zhang_UnderstandingClass_23}), and we show that they can be equivalently generated by minimizing}
a cost functional that integrates a weighted combination of state deviation and control effort over time \cite{liberzon_CalculusVariations_11, bertsekas_DynamicProgramming_00}, thereby achieving the best possible trade-off between tracking accuracy and control effort. By minimizing this functional, one can systematically shape the time evolution of the system to balance fast convergence with desirable transient properties.

In particular, we focus on a broad class of non-smooth, continuous-time, distributed algorithms \lic{that have been reported in the literature}
\cite{yamashita_PassivitybasedGeneralization_20}, which generalize augmented Lagrangian \lic{methods.}
We \lic{show} 
that these augmented dynamics can be generated by solving a corresponding network-wide optimal control problem, in which \li{a} 
cost functional penalizes deviations from the optimal solution during transients, \lhh{thereby serving as an exact performance metric for the particular \lic{algorithm.}} 
The augmentations introduce additional terms into the cost functional, which must be driven to zero by the optimal controller, thus providing an interpretation of the improved performance offered by these modifications.
\li{Variational principles have been used to provide interpretations} for accelerated gradient methods \cite{wibisono_VariationalPerspective_16} \lhh{and unaugmented network congestion control protocols \cite{lavaei_UtilityFunctionals_10}}.
\lic{Our work} is, to the best of our knowledge, the first to provide a network-wide performance guarantee for \li{general classes of} augmented primal–dual algorithms,
\li{
allowing nonlinear and non-smooth dynamics arising from general convex objective functions, and 
inequality \ilc{constraints.}
}

\lhk{This paper is organized as follows.}
We begin in Section \ref{sec_distributed_optimisation} by formulating \li{distributed optimization \lic{algorithms via appropriate equality} 
constraints on coupling variables}, and introduce the corresponding primal-dual algorithm and \lic{generalizations via augmentations} \cite{yamashita_PassivitybasedGeneralization_20}. In Section \ref{sec_opt_control_interp}, we present our main result, describing how the augmented algorithm can be generated by an optimal control problem via a meaningful cost functional. \lic{We start in Section \ref{sec_opt_control_interp} \lic{with a simplified problem formulation, and
}}Section \ref{sec_further_results} then \lic{presents} 
refinements and extensions to accommodate various problem settings found in the literature. Specifically,
\lic{we additionally introduce} \lic{coupling} 
inequality constraints, and extend the results to incorporate feed-forward terms.
Finally, \li{we 
conclude} in Section \ref{sec_conclusion}.
\lhk{
\ils{The proofs of the results} are given in the appendix.}

\subsubsection*{Notation and Definitions}
The set
$\mathbb{R}_{\geq 0}$ denotes the set of non-negative real numbers. The matrix-weighted \lhi{(semi-)}norm of a vector $x \in \mathbb{R}^n$ is given by $\lVert x \rVert_R^2 = x^T R x$ for positive \lhi{(}\li{semi}\lhi{)}definite $R \in \mathbb{R}^{n \times n}$.
Let $\mathcal{K}=\{\kappa_1, \ldots, \kappa_{|\mathcal{K}|}\}$ be an ordered index set. Then, the direct sum of a set of indexed matrices $B_k$ with $k \in \mathcal{K}$ is denoted by $\oplus_{k \in \mathcal{K}} B_k = \mathrm{diag} [B_{\kappa_1}, \ldots, B_{\kappa_{|\mathcal{K}|}}]$. In addition, the composite vector constructed from a set of indexed vectors $a_k$, where $k \in \mathcal{K}$, is denoted $[a_k]_{k \in \mathcal{K}} =[a_{\kappa_1}^T, \ldots,a_{\kappa_{|\mathcal{K}|}}^T]^T$.
$I_n$ represents the $n$-dimensional identity matrix, while $\mathbb{1}_n$ and $\mathbb{0}_n$ \lhn{represent} an $n$-dimensional vector of ones and zeros respectively.
\lhb{For a convex, continuously differentiable function $F:\Omega \rightarrow \mathbb{R}$, where $\Omega$ is a convex set, the Bregman divergence \cite{censor_ParallelOptimization_97,beck_MirrorDescent_03,banerjee_ClusteringBregman_05} associated with F between points $p, q \in \Omega$ is defined as
\begin{equation} \label{eq_bregman_def}
    D_F(p,q) = F(p) - F(q) - \nabla F(q)^T (p-q).
\end{equation}
We note that $D_F(p,q)$ is \li{non-negative} 
due to the convexity of $F(\cdot)$ \cite{boyd_ConvexOptimization_04}.}

\section{Augmented Distributed Optimization} \label{sec_distributed_optimisation}
\subsection{Static Optimization Problem}
\ic{We start by describing the 
\li{network} optimization problem under consideration.} \li{Further extensions to the formulation in this section will be described in the sections that follow.}
Consider a network \lhj{modeled} as a directed, connected graph $\mathcal{G}(\mathcal{V},\mathcal{E})$, where $\mathcal{V} = \{ \nu_1, \nu_2, \ldots, \nu_{|\calv|} \}$ is the set representing the nodes and $\mathcal{E} =  \{ \varepsilon_1, \varepsilon_2, \ldots, \varepsilon_{|\cale|} \}$ is the set representing the edges of the graph.
\lha{To each node $i \in \calv$, we assign an agent with an} objective to reach the minimum of a separable static cost function of the form
\begin{mini}
    {\tilde{\theta} \in \mathbb{R}^n}{\sum_{i \in \calv} F_i(\tilde{\theta}),}{\label{opt_distributed}}{}
        \addConstraint {G_i(\tilde{\theta}) & \leq 0, \qquad i \in \calv, }
\end{mini}
where the functions $F_i:\mathbb{R}^n \rightarrow \mathbb{R}$ are convex, continuously differentiable objective functions associated with the agent $i \in \calv$ \cite{yang_SurveyDistributed_19}.
The functions $G_i = [G_{i1}, G_{i2}, \ldots, G_{im}]^T :\mathbb{R}^{n} \rightarrow \mathbb{R}^m$ (with each $G_{ik}:\mathbb{R}^{n} \rightarrow \mathbb{R} , k = \{1, \ldots, m\}$ being convex and continuously differentiable) correspond to local
inequality conditions on the global optimization variable $\tilde{\theta} \in \mathbb{R}^n$ at each \li{node.}
\lhab{
\lic{To simplify the presentation, we assume in what follows that} $n = m =1$ so that $\tilde{\theta} \in \mathbb{R}$.}

In order for each agent to arrive at the same optimal value for the global variable $\tilde{\theta}^{\ast} \lhab{\in \mathbb{R}}$ in a distributed fashion, some type of communication and consensus protocol is required between agents. Let the links in this communication network be represented by the set $\cale$. In a distributed algorithm, each agent $i \in \calv$ updates its own local estimate of the solution $\theta_i \in \mathbb{R}$ based on local information.
Such algorithms can be constructed via the following transformation of \licc{the} optimization problem \eqref{opt_distributed}:
\begin{mini}
	{\theta \in \mathbb{R}^{|\calv|}}{\sum_{i \in \calv} F_i(\theta_i),}{\label{opt_distributed_transformed}}{}
        \addConstraint{G_i(\theta_i) & \leq 0, \qquad i \in \calv }
	\addConstraint { \cala^T \theta & = 0,}
\end{mini}
where $\theta = [\theta_i]_{i \in \calv}$ and $\cala \in \mathbb{R}^{|\calv| \times |\cale|}$ is the matrix
\lic{that
enforces the consensus constraint $\theta_p = \theta_q$ for all $p, q \in \calv$. The sparsity structure of $\cala$ also determines the communication structure of the distributed algorithm that will be constructed  (e.g., $\cala$ can be a Laplacian matrix, or an incidence matrix \cite{droge_ContinuoustimeProportionalintegral_14}).}

\begin{remark}
    \lhi{}
    \icl{We would like to note that for simplicity in the presentation, we have considered \ils{in \eqref{opt_distributed_transformed}} 
    a simplified version of the optimization problem}
    \lhi{with one coupling variable, i.e., $n=1$ in \eqref{opt_distributed}.
    \lic{In the case $n>1$, i.e., multiple coupling variables are present, a formulation analogous to \eqref{opt_distributed_transformed} can be constructed. That is, each agent has separate variables $\theta_i\in\mathbb{R}^{n_i}$ with a corresponding equality constraint introduced of the form $\tilde \cala^T \theta = 0$, $\theta = [\theta_i]_{i\in\calv}$ to ensure equivalence with \eqref{opt_distributed}. }}
\end{remark}
\begin{remark}
    \lhb{\li{Instead} of a consensus constraint in \eqref{opt_distributed_transformed}, it is possible to modify the optimization problem to \lic{include} \lhi{coupling} 
    inequality \li{constraints} of the \lic{form $H(\theta) \leq 0$, where $H(\cdot)$
    is a convex continuously differentiable function}
    (such constraints appear, for example, in \lic{congestion control problems in} 
    communication networks \cite{kelly_RateControl_98}). The results \lic{described in this section} 
    can be easily extended to include problems of this form, \lhe{which is} 
    \lic{described} in Section \ref{sec_global_ineq}.}
\end{remark}

\lhh{The optimization problem \eqref{opt_distributed_transformed} admits the following Lagrangian function
\begin{equation} \label{eq_lagrangian}
    L_0(\theta,\lambda,\mu) = F(\theta) + \lambda^T G(\theta) + \mu^T \cala^T \theta,
\end{equation}
where $F(\theta) = \sum_{i \in \calv} F_i(\theta_i)$,
\lhj{$G(\theta) = [G_i(\theta_i)]_{i\in\calv}$,}
$\lambda := [\lambda_i]_{i \in \calv} \in \mathbb{R}^{|\calv|}$
is the dual variable (or Lagrange multiplier) associated with the inequality constraints, and $\mu := [\mu_j]_{j \in \cale} \in \mathbb{R}^{|\cale|}$ is the dual variable associated with the consensus constraint. }
\lic{We also assume Slater's condition holds\lhj{\cite{boyd_ConvexOptimization_04}} and 
\lhj{$F(\theta)$} \lic{is bounded from below} 
    in the feasible set of \eqref{opt_distributed_transformed}.}


\lhi{\lic{Saddle points $(\theta^\ast, \lambda^\ast, \mu^\ast)$ of \eqref{eq_lagrangian} are then optimal solutions to the optimization problem \eqref{opt_distributed_transformed}.   \lic
In particular, these satisfy} 
the}
\lhb{Karush–Kuhn–Tucker (KKT) \li{optimality} conditions \cite{boyd_ConvexOptimization_04}
}
\begin{subequations} \label{eq_distributed_kkt}
    \begin{align}
        & \nabla F(\theta^{\ast}) + \nabla G(\theta^{\ast}) \lambda^{\ast} + \cala \mu^{\ast} = 0 \label{eq_distributed_kkt_1} \\
        & \lambda_i^{\ast}G_i(\theta_i^{\ast}) = 0, \quad \lambda_i^{\ast} \geq 0, \quad G_i(\theta_i^{\ast}) \leq 0, \quad \forall i \in \calv \label{eq_distributed_kkt_2} \\
        & \cala^T \theta^{\ast} = 0, \label{eq_distributed_kkt_3}
    \end{align}
\end{subequations}
\lhj{where $\nabla G(\theta^\ast) = \mathrm{diag}[\nabla G_i(\theta_i^\ast)]_{i\in\calv}$ is the Jacobian of $G(\cdot)$ evaluated at $\theta^\ast$.}

\subsection{Primal-Dual Algorithm}
\lha{For strictly convex cost functions,} this problem can be solved in a distributed fashion using the primal-dual algorithm \cite{yang_SurveyDistributed_19}, which converges \ic{to an optimal solution $(\theta^{\ast}, \lambda^{\ast}, \mu^{\ast})$
via} the following continuous dynamics for the primal variable $\theta$ and the dual variables $\lambda$ and $\mu$:
\begin{subequations} \label{eq_primal_dual_alg}
    \begin{align}
        \dot{\theta} &= -\nabla F(\theta) - \nabla G(\theta) \lambda - \cala \mu, \label{eq_primal_dual_alg_x} \\
        \dot{\lambda} & = [G(\theta)]_{\lambda}^+, \ \lambda(0) \geq 0, \label{eq_primal_dual_alg_lambda} \\
        \dot{\mu} &= \cala^T \theta. \label{eq_primal_dual_alg_mu}
    \end{align}
\end{subequations}
where the $k^{\mathrm{th}}$ component of the operator $[\cdot]_{\ast}^+$ is defined by
\begin{equation} \label{eq_positive_projection}
    \left( \left[ \sigma \right]^+_{\epsilon} \right)_k := \begin{cases}
        \sigma_k, & \textrm{if $\epsilon_k>0$,} \\
        \sigma_k, & \textrm{if $\epsilon_k = 0$ and $\sigma_k \geq 0$,} \\
        0, & \textrm{if $\epsilon_k = 0$ and $\sigma_k < 0$.}
    \end{cases}
\end{equation}
for vectors $\sigma$ and $\epsilon$ of equal dimension.
Here, we assume $\nabla F(\cdot)$ and all $\nabla G_{ik}(\cdot)$ are locally Lipschitz.
\lha{We see that the KKT conditions for optimality \eqref{eq_distributed_kkt} are satisfied by the equilibrium points of the algorithm \eqref{eq_primal_dual_alg}.}

\begin{figure}[t]
\centering
\small
\input{Images/tikz_block_diagram}
\caption{Closed loop block diagram of the augmented primal-dual distributed optimization algorithm \eqref{eq_primal_dual_alg}, \lhb{with primal dynamics given by block-diagonal system $\Sigma_\theta$, dual dynamics associated with local inequality constraints given by block-diagonal $\Sigma_\lambda^+$, and dual dynamics associated with the consensus constraint given by block-diagonal $\Sigma_\mu$}.}
\label{fig_dist_opt}
\end{figure}
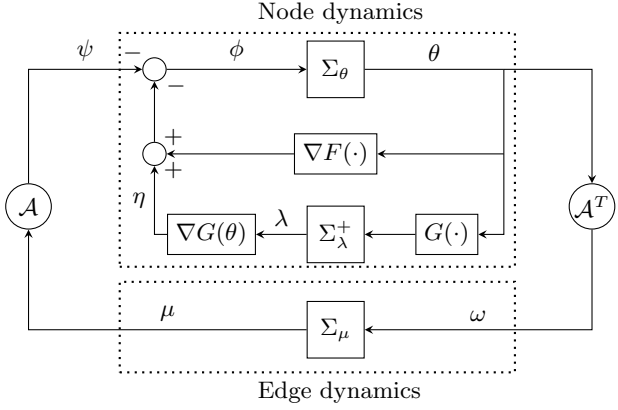

\icl{\iccc{Now, the system \eqref{eq_primal_dual_alg} can be reformulated as an interconnection of \lha{dynamical systems associated with each} node and edge,
as described below and illustrated in Figure \ref{fig_dist_opt}. 
We define the variables}
\begin{subequations} \label{eq_var_defs}
    \begin{align}
        \omega &= [\omega_j]_{j \in \cale} := \cala^T \theta \label{eq_def_z}\\
        \eta&= [\eta_i]_{i \in \calv} := \nabla G(\theta) \lambda \label{eq:def_eta}\\
        \psi&= [\psi_i]_{i \in \calv} :=  \cala \mu \label{eq_def_psi}\\
        \phi &= [\phi_i]_{i \in \calv} := - \nabla F(\theta) - \eta - \psi \label{eq_def_phi}
    \end{align}
\end{subequations}
and their corresponding optimal values $\omega^{\ast} := [\omega_j^{\ast}]_{j \in \cale} =  \cala^T \theta^{\ast}$ (so we have $\omega^{\ast} = 0$ via \eqref{eq_distributed_kkt_3}), $\eta^{\ast}:= [\eta_i^{\ast}]_{i \in \calv} =  \nabla G(\theta^{\ast}) \lambda^{\ast} $,  $\psi^{\ast} := [\psi_i^{\ast}]_{i \in \calv}= \cala \mu^{\ast}$ \lhg{and $\phi^{\ast} := \lhf{[\phi_i^{\lhad{\ast}}]_{i \in \calv}} = \nabla F(\theta^{\ast}) + \eta^{\ast} + \psi^{\ast}$ (so we have that $\phi^{\ast} = 0$ via \eqref{eq_distributed_kkt_1}).}}



\icl{We consider here as node subsystems the systems with input \lhg{$\psi_i$} and output \lhg{$\theta_i$}, as illustrated in Figure \ref{fig_dist_opt}. In particular, these are described by the interconnection of the primal dynamics \eqref{eq_primal_dual_alg_x} with the dual inequality constraint dynamics \eqref{eq_primal_dual_alg_lambda}. We consider as edge subsystems the systems in Figure \ref{fig_dist_opt} with input \lhg{$\omega_i$} and output \lhg{$\mu_i$}, and these are described by \eqref{eq_primal_dual_alg_mu}.}

\icl{\icl{In Figure \ref{fig_dist_opt}, the \lhab{operators}}
\lhm{$\Sigma_\theta := \bigoplus_{i \in \calv} \Sigma_{\theta i}$, $\Sigma_\lambda^+:= \bigoplus_{i \in \calv} \Sigma_{\lambda i}^+$ and $\Sigma_\mu:= \bigoplus_{j \in \cale} \Sigma_{\mu j}$ 
have inputs and outputs as illustrated in the figure 
with these related via
the dynamics in \eqref{eq_primal_dual_alg}. In particular, \lhab{to recover \eqref{eq_primal_dual_alg}, we see that each} $\Sigma_{\theta i}$,  $\Sigma_{\mu j}$ are integrators and $\Sigma_{\lambda i}^+$ are nonlinear operators given \icc{by 
\eqref{eq_primal_dual_alg_lambda}.}}}
\subsection{Augmented Algorithm} \label{sec_augmented_algorithm}
\lha{Using an argument based on passivity theory \cite{khalil_NonlinearSystems_13}, 
\lho{it was shown in \cite{yamashita_PassivitybasedGeneralization_20}}
that the algorithm \eqref{eq_primal_dual_alg} can be augmented with auxiliary dynamics at each node and edge. With these augmentations, the primal-dual algorithm can be used to solve distributed optimization problems of the form \eqref{opt_distributed_transformed} with non-strictly convex cost functions.}
\lhm{\lhab{As in \cite{yamashita_PassivitybasedGeneralization_20}, we} introduce auxiliary dynamics \lhg{for the systems associated with} each $i \in \calv$ and each $j \in \cale$ \lha{so that the operators $\Sigma_{\theta i}$, $\Sigma_{\lambda i}^+$ and $\Sigma_{\mu j}$ in Figure \ref{fig_dist_opt} are} as follows:
\begin{subequations}\label{eq_comp_dynamics}
    \begin{align}
        & \Sigma_{\theta i} : \begin{cases}
            \dot{\xi}_{i1} &= b_{i1} \phi_i, \\
            \dot{\xi}_{ik} &= b_{ik} \phi_i - a_{ik} \xi_{ik}, \quad k = 2, \ldots, \rho_i,
        \end{cases} \label{eq_comp_dynamics_xi} \\
        & \Sigma_{\lambda i}^+ : \begin{cases}
            \dot{\tau}_{i1} &= [\hat{b}_{i1} G_i(\theta_i)]_{\tau_{i1}}^+, \\
            \dot{\tau}_{ik} &= [\hat{b}_{ik} G_i(\theta_i) - \hat{a}_{ik} \tau_{ik}]_{\tau_{ik}}^+, \quad k = 2, \ldots, \hat{\rho}_i,
        \end{cases} \label{eq_comp_dynamics_tau} \\
        & \Sigma_{\mu j} : \begin{cases}
            \dot{\zeta}_{j1} &= \check{b}_{j1} \omega_j, \\
            \dot{\zeta}_{jk} &= \check{b}_{jk} \omega_j - \check{a}_{jk} \zeta_{jk}, \quad k = 2, \ldots, \check{\rho}_{\lhn{j}},
        \end{cases} \label{eq_comp_dynamics_zeta}
    \end{align}
\end{subequations}
where $\ils{\xi_{ik}(t)} \in \mathbb{R}, k \in \{1, \ldots, \rho_i \}$, $\ils{\tau_{ik}(t)} \in \mathbb{R}_{\geq 0}, k \in \{1, \ldots, \hat{\rho}_i \}$ and $\ils{\zeta_{jk}(t)} \in \mathbb{R}, k \in \{1, \ldots, \check{\rho}_{\lhn{j}} \}$ are \textit{auxiliary states} and \icc{the outputs of systems $\Sigma_{\theta i}, \Sigma_{\lambda i}^+, \Sigma_{\mu j}$ are, respectively,}
\begin{subequations} \label{eq_comp_outputs}
\begin{align}
    &\theta_i = \textstyle \sum_{k = 1}^{\rho_i} \xi_{ik} + d_i \phi_i, \label{eq_comp_outputs_x} \\
    & \lambda_i = \textstyle \sum_{k = 1}^{\hat{\rho}_i} \tau_{ik} + \hat{d}_i \mathrm{max}\{ 0, G_i(\theta_i)\}, \label{eq_comp_outputs_lambda} \\
    & \mu_j =\textstyle \sum_{k = 1}^{\check{\rho}_j} \zeta_{jk} + \check{d}_j \omega_j. \label{eq_comp_outputs_mu}
\end{align}
\end{subequations}
Here,} the constants satisfy $a_{i\rho_i} > \cdots > a_{i2} > 0$, $b_{ik} > 0, \forall k \in \{1, \ldots, \rho_i\}$, and $d_i \geq 0$, \lhab{with similar conditions holding for $\hat{a}_{ik}$, $\hat{b}_{ik}$, $\hat{d}_i$, $\check{a}_{jk}$, $\check{b}_{jk}$ and $\check{d}_j$.}
\lha{The terms containing $d_i$, $\hat{d}_i$ and $\check{d}_j$ in \eqref{eq_comp_outputs} are termed \textit{feed-forward} terms, as they directly add the respective inputs of the subsystems in Figure \ref{fig_dist_opt} to their outputs (multiplied by a \lhe{non-negative} constant, and with corresponding max operator applied in the case of $\lambda_i$).}


\icl{Note that the dynamics in \lhg{\eqref{eq_comp_dynamics}}, \eqref{eq_comp_outputs} are \lhab{linear} combinations of integrators and phase lead compensators \lhab{(see frequency domain description in \cite{yamashita_PassivitybasedGeneralization_20})}, used in conjunction with appropriate projections in \eqref{eq_comp_dynamics_tau}, \eqref{eq_comp_outputs_lambda} to ensure positivity \icc{of $\tau_i$, $\lambda_i$.}}

\begin{remark} \label{remark_alg_equilibrium}
    \lhc{We note that the algorithm \eqref{eq_primal_dual_alg} is only guaranteed \icl{to 
    converge} to the optimal solution \icl{of the 
    optimization problem} \eqref{opt_distributed_transformed} if $F(\cdot)$ is strictly convex. By contrast, the augmented algorithm with lead compensators \eqref{eq_comp_dynamics} \lhm{and \eqref{eq_comp_outputs}} \icl{
    converges} to a particular optimal solution $(\theta^{\ast},\lambda^{\ast},\mu^{\ast})$ (\icl{that depends} on the algorithm initial conditions) even in the case where $F(\cdot)$ is convex but not strictly convex, provided $\rho_i \geq 2$ or $d_i > 0$
    for all $i \in \calv$ \cite{yamashita_PassivitybasedGeneralization_20}.}
\end{remark}
\subsection{Improved Transient Performance} \label{sec_improved_performance}

\lha{We observe that the modified algorithm \eqref{eq_comp_dynamics} reduces to the standard primal-dual algorithm \eqref{eq_primal_dual_alg} in the case where each $\rho_i$, $\hat{\rho}_i$, $\check{\rho}_j$ is equal to 1, and each $d_i$, $\hat{d}_i$, $\check{d}_i$ is equal to 0, for all $i \in \calv$ and $j \in \cale$. 
\lho{It was shown in \cite{yamashita_PassivitybasedGeneralization_20}}
that the auxiliary states introduced in \eqref{eq_comp_dynamics} and the feed-forward terms in \eqref{eq_comp_outputs} can be viewed as generalizations of augmented Lagrangian methods \lhi{(see examples section below)},
which have been demonstrated via simulation to improve the performance of the standard primal-dual algorithm \eqref{eq_primal_dual_alg}. However, a theoretical performance metric to quantify these improvements has been lacking in the literature.}
\lha{We now use a simple example to exhibit the improved performance offered by these algorithms.}

\subsubsection*{Example}
\begin{figure}[t]
    \centering
    \input{Images/tikz_network}
    \caption{Graph structure of the example used in Section \ref{sec_improved_performance}.}
    \label{fig_network}
\end{figure}
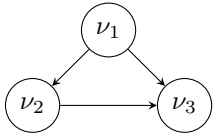

\lha{We consider a simple example using the network shown in Figure \ref{fig_network}, which consists of three nodes and three edges. We let the cost functions equal:
\begin{subequations}
    \begin{align}
        F_{\nu 1}(\theta_{\nu 1}) & = (\theta_{\nu 1} - 0.5)^2, \\
        F_{\nu 2}(\theta_{\nu 2}) & = e^{-0.5\theta_{\nu 2}}, \\
        F_{\nu 3}(\theta_{\nu 3}) & = - \log{\theta_{\nu 3}},
    \end{align}
\end{subequations}
and we ignore the local inequality constraints so that the system $\Sigma_\lambda^+ = 0$. The matrix $\cala$ is taken as the directed network incidence matrix (where each node is given arbitrary direction), with entry $(i,j)$ equal to $1$ when $i \in \calv$ is the sink node for edge $j \in \cale$, $-1$ when $i \in \calv$ is the source node for edge $j \in \cale$, and $0$ otherwise.}

\lha{We test two different augmentations of the algorithm \eqref{eq_comp_dynamics}:}
\begin{enumerate}
    \item \lha{\textbf{Auxiliary states:} For each $i \in \calv$, we let $\rho_i=2$ and for each $j\in \cale$, we let $\check{\rho}_j=1$, so that the system associated with each node has one \lhb{additional} auxiliary state contributing to the value of $\theta_i$ via \eqref{eq_comp_outputs_x}. We let $d_i = \check{d}_j=0$ so that there are no input feed-\lhj{forward} terms in the outputs \eqref{eq_comp_outputs_x} and \eqref{eq_comp_outputs_mu}. The parameters $b_{i1}$ and $b_{i2}$ \lhh{are set to $\frac{1}{2}$}, $a_{i2}$ \lhh{set to 2} and $\check{b}_{j1}$ \lhh{set to 1}.}
    \lhh{This gives the following algorithm:
        \begin{align*}
            &\Sigma_{\theta i} : \begin{cases}
                \dot{\xi}_{i1} &= \frac{1}{2}\phi_i, \\
                \dot{\xi}_{i2} &= \frac{1}{2}\phi_i - 2 \xi_{i2}, \\
                \theta_i & = \xi_{i1} + \xi_{i2}, \\
            \end{cases} & i \in \{\nu_1,\nu_2,\nu_3\},  \\
            &\Sigma_{\mu_j} : \begin{cases}
                \dot{\mu}_{j} &= \omega_j, \\
            \end{cases} & j \in \{ \varepsilon_1,\varepsilon_2,\varepsilon_3\},
        \end{align*}
    where $\omega_i$ and $\phi_j$ are given by the network interconnection, as defined in \eqref{eq_def_z} and \eqref{eq_def_phi} respectively.
    By defining a new variable $\theta_i^\prime = \xi_{i1}-\xi_{i2}$ for each $i \in \calv$, we obtain the following equivalent algorithm:
        \begin{align*}
            &\Sigma_{\theta i} : \begin{cases}
                \dot{\theta}_i &= \phi_i + (\theta_i^\prime - \theta_i), \\
                \dot{\theta}_i^\prime &= \theta_i-\theta_i^\prime,
            \end{cases} & i \in \{\nu_1,\nu_2,\nu_3\},  \\
            &\Sigma_{\mu_j} : \begin{cases}
                \dot{\mu}_{j} &= \omega_j, \\
            \end{cases} & j \in \{ \varepsilon_1,\varepsilon_2,\varepsilon_3\}.
        \end{align*}
    Compared to \eqref{eq_primal_dual_alg}, this transformed algorithm contains an additional auxiliary state at each node $i \in \calv$. These dynamics can be obtained by deriving the primal-dual dynamics that seek the saddle point of the augmented Lagrangian
    \begin{equation}
        L_a(\theta,\theta^\prime,\mu) = F(\theta) + \textstyle \frac{1}{2} \lVert \theta^\prime - \theta \rVert^2 + \mu^T \cala^T \theta,
    \end{equation}
    where $\theta^\prime = [\theta_{\nu1}^\prime,\theta_{\nu2}^\prime,\theta_{\nu3}^\prime]^T$. This is the augmented Lagrangian used in proximal point approaches \cite{gill_PrimaldualAugmented_12, rockafellar_AugmentedLagrangians_76, holding_StabilityInstability_15}.}

    \item \lha{\textbf{Feed-forward term:} For each $i \in \calv$ and $j \in \cale$, we let $\rho_i$, $\check{\rho}_j$ equal to 1 (so there are no \lhb{additional} auxiliary states) and take $b_{i1}=\check{b}_{j1}=1$. In \eqref{eq_comp_outputs}, we let each $d_i = 0$ and each $\check{d}_j=1$, so the output of $\Sigma_\mu$ feeds the input $\omega$ forward in Figure \ref{fig_dist_opt}. As we will see in Section \ref{sec_feedforward} (see Figure \ref{fig_transformation}), this is equivalent to using \lhh{Laplacian} feedback of $\theta$ in the system $\Sigma_\theta$, giving rise to the following augmentation of \eqref{eq_primal_dual_alg}:
    \begin{subequations} \label{eq_primal_dual_lagrangian}
        \begin{align}
            \dot{\theta} &= -\nabla F(\theta) - \cala \mu - \mathcal{L} \theta, \\
            \dot{\mu} &= \cala^T \theta,
        \end{align}
    \end{subequations}
    where $\mathcal{L} = \cala\cala^T$ is the network \lhe{Laplacian}.}
    \lhh{This algorithm (sometimes known as the PI distributed optimization algorithm \cite{droge_ContinuoustimeProportionalintegral_14}) can equivalently be obtained by deriving the primal-dual dynamics that seek the saddle point of the augmented Lagrangian
    \begin{equation}
        L_f(\theta,\mu) = F(\theta) + \textstyle \frac{1}{2} \lVert \cala^T \theta \rVert^2 + \mu^T \cala^T \theta,
    \end{equation}
    which is an application of the Hestenes-Powell approach \cite{hestenes_MultiplierGradient_69, powell_MethodNonlinear_69, rockafellar_AugmentedLagrange_74, bertsekas_ConstrainedOptimization_82} to a distributed optimization problem of the form \licc{in} \eqref{opt_distributed_transformed}.}
\end{enumerate}

\begin{figure}[t]
\centering
\includegraphics[width=0.45\textwidth]{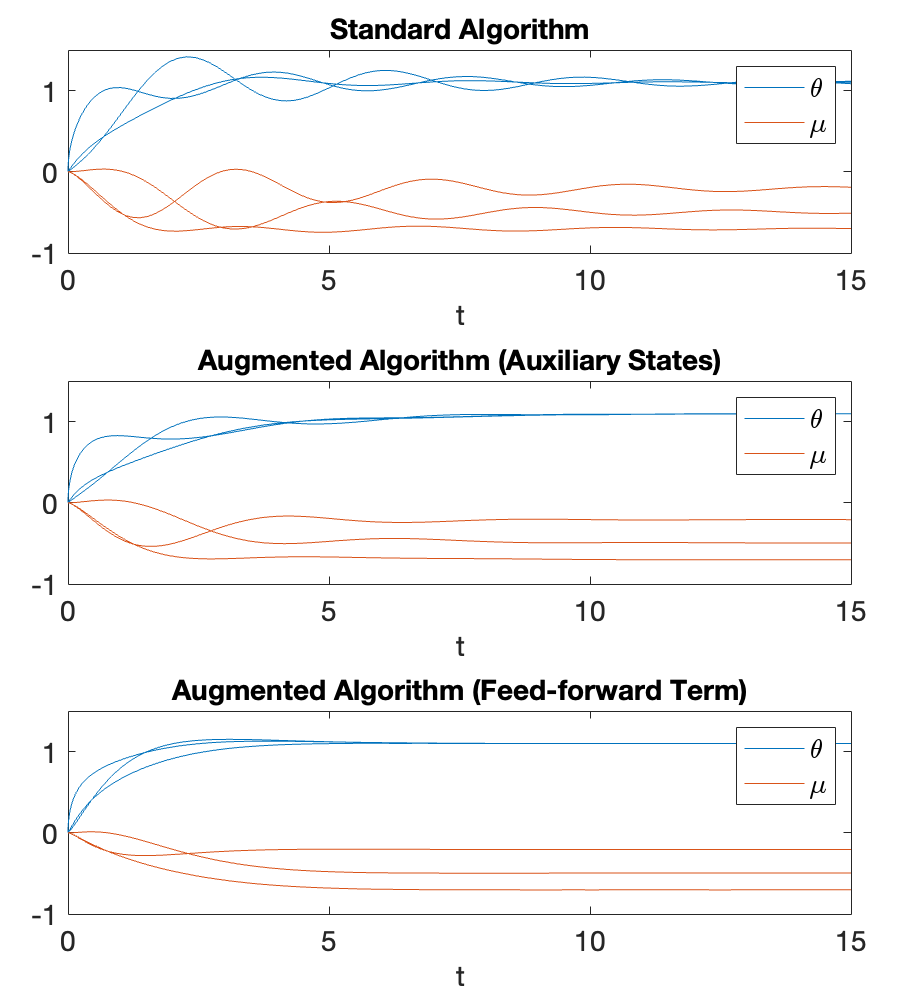}
\caption{Case study demonstrating the improved performance offered by using two different augmentations of the primal-dual algorithm, which are described by the generalized dynamics \eqref{eq_comp_dynamics}. The standard algorithm is given by \eqref{eq_primal_dual_alg}; the auxiliary states case adds an additional state contributing to the output $\theta$ at each node; and the feed-forward term case is \lic{equivalent} to using Lagrangian feedback in the primal dynamics. }
\label{fig_example}
\end{figure}

\lha{The results of using the standard algorithm \licc{in} \eqref{eq_primal_dual_alg} and the two augmentations described above are shown in Figure \ref{fig_example}. We can see that using these simple augmentations greatly improves the transient performance of the algorithm by damping oscillations and improving algorithm \li{settling time}.}

\section{Optimal Control Interpretation} \label{sec_opt_control_interp}
\lha{We now show that the performance \li{improvements}  
offered by the augmented algorithm \eqref{eq_comp_dynamics} can be quantified via an optimal control interpretation \cite{liberzon_CalculusVariations_11,bertsekas_DynamicProgramming_00}. \lhe{As discussed in the Introduction, the goal of optimal control is to determine the best possible control input for a dynamical system that minimizes a cost functional consisting of a weighted combination of state deviation and control effort over the transient response of the dynamical system. A typical infinite horizon optimal control problem takes the form
\begin{equation}
    \begin{aligned}
        \min_{u} \quad \int_0^\infty \lVert \ils{u(t)} \rVert_R^2 + q(\ils{x(t)}) \ dt, \\
        \text{s.t.} \qquad \ils{\dot{x} = f(x,u),}
    \end{aligned}
\end{equation}
where \ils{$x(t) \in \mathbb{R}^n$ is the system state, $u(t)\in \mathbb{R}^m$} 
is an \ils{arbitrary 
control} input, $R \in \mathbb{R}^{m \times m} $ is a positive definite, symmetric control cost weighting matrix, $q:\mathbb{R}^n \rightarrow \mathbb{R}$ is a positive semi-definite state cost function, and $f:\mathbb{R}^n \times \mathbb{R}^m \rightarrow \mathbb{R}^n$ is \ils{locally} Lipschitz continuous. We will show that the augmented algorithm \licc{in} \eqref{eq_comp_dynamics} can be generated by a problem of this form.}}

\lha{In order to simplify the presentation in what follows, we first consider the case with no feed-forward terms in the system outputs \eqref{eq_comp_outputs} (i.e., $d_i=0$, $\hat{d}_i=0$ and $\check{d}_j=0$ for all $i \in \calv$ and $j \in \cale$). The effect of the feed-forward terms will be presented in Section \ref{sec_feedforward}.}

\vspace{.25cm}
\subsection{Control Formulation}
Firstly, we need to form dynamics \lha{for the systems in Figure \ref{fig_dist_opt}} with a control input.
For each $i \in \calv$ and $j \in \cale$,
\lha{we \ils{formulate} 
dynamics in vector form as follows:}
\begin{subequations} \label{eq_primal_dual_control}
    \begin{align}
        &\Sigma_{\theta i} = \begin{cases}
            \dot{\xi}_i = b_i \mathbb{1}_{\rho_i} \phi_i + b_{ui} u_i, \\
            \theta_i = \mathbb{1}_{\rho_i}^T \xi_i,
        \end{cases} \label{eq_primal_dual_control_theta} \\
        &\Sigma_{\lambda i}^+ = \begin{cases}
            \dot{\tau}_i = [ \hat{b}_i \mathbb{1}_{\hat{\rho}_i} G_i(\theta_i) + \hat{b}_{ui} \hat{u}_i ]^+_{\tau_i}, \\
            \lambda_i = \mathbb{1}_{\hat{\rho}_i}^T \tau_i,
        \end{cases} \label{eq_primal_dual_control_lambda} \\
        &\Sigma_{\mu j} = \begin{cases}
            \dot{\zeta}_j = \check{b}_j \mathbb{1}_{\check{\rho}_j} \omega_j + \check{b}_{uj} \check{u}_j, \\
            \mu_j = \mathbb{1}_{\check{\rho}_j}^T \zeta_j.
        \end{cases} \label{eq_primal_dual_control_mu}
    \end{align}
\end{subequations}
\lha{Here, we define:
state vectors} $\xi_i := [\xi_{i1}, \ldots, \xi_{i \rho_i}]^T$, $\tau_i := [\tau_{i1}, \ldots, \tau_{i \hat{\rho}_i}]^T$, and $\zeta_j := [\zeta_{j1}, \ldots, \zeta_{j \check{\rho}_j}]^T$;
\lha{state matrices} $b_i := \mathrm{diag} [b_{i1},\ldots, b_{i\rho_i}]$ (with $b_{ik} > 0, \forall k \in \{1, \ldots, \rho_i\}$), $\hat{b}_i := \mathrm{diag} [\hat{b}_{i1},\ldots, \hat{b}_{i\hat{\rho}_i}]$ (with $\hat{b}_{ik} > 0, \forall k \in \{1, \ldots, \hat{\rho}_i\}$), and $\check{b}_j := \mathrm{diag} [\check{b}_{j1},\ldots, \check{b}_{j \check{\rho}_j}]$ (with $\check{b}_{jk} > 0, \forall k \in \{1, \ldots, \check{\rho}_j\}$);
\lha{and control input matrices} $b_{ui} = I_{o \rho_i}$, $\hat{b}_{ui} = I_{o \hat{\rho}_i}$, and $\check{b}_{uj} = I_{o \check{\rho}_j}$, where
$I_{o n}:= \begin{bmatrix} \mathbb{0}_{n-1} & I_{n - 1} \end{bmatrix}^T$ for $n \in \mathbb{N}$.
\lha{In addition,}
\lhg{$\omega_j$ and $\phi_i$ \li{are as} defined in \eqref{eq_def_z} and \eqref{eq_def_phi} respectively.}

We \ils{let  $u_i(t)\in \mathbb{R}^{\rho_i - 1}$, $\hat{u}_i(t)\in \mathbb{R}^{\hat{\rho}_i - 1}$ and $\check{u}_j(t)\in \mathbb{R}^{\check{\rho}_j-1}$} be 
\li{
\lic{
\ils{control inputs that are stabilizing, i.e. the system states converge to an equilibrium point where}
\licc{variables $(\theta,\lambda,\mu$)} 
are equal to the optimal solution \lhj{$(\theta^\ast,\lambda^\ast,\mu^\ast)$, which is the same solution} the primal-dual dynamics in \eqref{eq_comp_dynamics} converge to}. We}
observe that we recover the augmented primal-dual dynamics \eqref{eq_comp_dynamics} if we let $u_i$, $\hat{u}_i$ and $\check{u}_j$ equal to
\begin{subequations} \label{eq_augmented_controllers}
    \begin{align}
        \underline{u}_i &= \left[ -a_{i2} \xi_{i2}, \ldots, -a_{i\rho_{i}} \xi_{i\rho_{i}} \right]^T,&    \forall i \in \calv, \label{eq_augmented_controllers_theta} \\
        \underline{\hat{u}}_i &= \left[ -\hat{a}_{i2} \tau_{i2}, \ldots, -\hat{a}_{i\hat{\rho}_{i}} \tau_{i\hat{\rho}_{i}} \right]^T,&    \forall i \in \calv, \label{eq_augmented_controllers_lambda}  \\
        \underline{\check{u}}_j &= \left[ -\check{a}_{j2} \zeta_{j2}, \ldots, -\check{a}_{j\check{\rho}_j} \zeta_{j\check{\rho}_j} \right]^T, &  \forall j \in \cale. \label{eq_augmented_controllers_mu}
    \end{align}
\end{subequations}

\lha{We now state \lho{one of the main results} of the paper, which is proved in Appendix \ref{sec_proof_p1}.
}

\lha{\begin{proposition} \label{prop_main_result}
  Consider the following optimal control problem for system \eqref{eq_primal_dual_control} with initial condition $(\xi_0, \tau_0, \zeta_0)$:
    \begin{equation} \label{eq_primal_cost}
        \begin{split}
            \min_{U} \int_0^{\infty} & \lVert U \rVert_R^2 + \sum_{i \in \calv} q_i(\xi_i,\tau_i, R_i)   + \sum_{i \in \calv} \hat{q}_i(\tau_i,\hat{R}_i) \\& + \sum_{j \in \cale} \check{q}_j(\zeta_j,\check{R}_j) \ dt,
        \end{split}
    \end{equation}
    \ils{where $U = \begin{bmatrix} [u_i]_{i \in \calv}^T & [\hat{u}_i]_{i \in \calv}^T & [\check{u}_j]_{j \in \cale}^T \end{bmatrix}^T$. 
    Then} $\underline{U} = \begin{bmatrix} [\underline{u}_i]_{i \in \calv}^T & \underline{\hat{u}}_i]_{i \in \calv}^T & [\underline{\check{u}}_j]_{j \in \cale}^T \end{bmatrix}^T$ given by \eqref{eq_augmented_controllers} is the optimal solution to the optimal control problem \eqref{eq_primal_cost}, with $R= \mathrm{diag} \begin{bmatrix} \oplus_{i \in \calv} R_i & \oplus_{i \in \calv} \hat{R}_i & \oplus_{j \in \cale} \check{R}_j \end{bmatrix}$,
    where
    \begin{subequations} \label{eq_primal_dual_R}
        \begin{align}
            & R_i = \textstyle \frac{1}{2} \mathrm{diag} \left[ \frac{1}{a_{i2}b_{i2}}, \ldots, \frac{1}{a_{i\rho_i}b_{i\rho_i}} \right], \\
            & \hat{R}_i = \textstyle  \frac{1}{2} \mathrm{diag} \left[ \frac{1}{\hat{a}_{i2}\hat{b}_{i2}}, \ldots, \frac{1}{\hat{a}_{i\hat{\rho}_i}\hat{b}_{i\hat{\rho}_i}} \right], \\
            & \check{R}_j = \textstyle  \frac{1}{2} \mathrm{diag} \left[ \frac{1}{\check{a}_{j2}\check{b}_{j2}}, \ldots, \frac{1}{\check{a}_{j\check{\rho}_j}\check{b}_{j\check{\rho}_j}} \right],
        \end{align}
    \end{subequations}
    and $q_i(\xi_i,\tau_i,R_i)$, $\hat{q}_i(\tau_i,\xi_i,\hat{R}_i)$ and $\check{q}_j(\zeta_j ,\check{R}_j)$ are given by~\eqref{eq_state_cost_primal_dual}.
\begin{figure*}[t]
    \begin{equation} \label{eq_state_cost_primal_dual}
        \begin{aligned}
                q_i(\xi_i,\tau_i,R_i)  & =  {\underbrace{%
                    \vphantom{ \textstyle \frac{a_{ik}}{2b_{ik}}} 
                    \left(\theta_i - \theta_i^{\ast} \right) \left(\nabla F_i(\theta_i) - \nabla F_i(\theta_i^{\ast}) \right)}_{\sigma_{i1}}
                }
                + {\underbrace{%
                    \vphantom{ \textstyle \frac{a_{ik}}{2b_{ik}}} 
                    \left(\theta_i - \theta_i^{\ast} \right) \left( \eta_i - \eta_i^{\ast} \right)}_{\sigma_{i2}}}
                + {\underbrace{%
                    \textstyle \sum_{k = 2}^{\rho_i} \frac{a_{ik}}{2b_{ik}} \xi_{ik}^2}_{\sigma_{i3}}},
            \\
                \hat{q}_i(\tau_i,\xi_i,\hat{R}_i)  &=
                {\underbrace{%
                    \vphantom{ \textstyle \frac{\hat{a}_{ik}}{2\hat{b}_{ik}}} 
                    -(\tau_{i1} - \lambda_i^{\ast})[G_i(\theta_i)]_{\tau_{i1}}^+
                    - \textstyle \sum_{k = 2}^{\hat{\rho}_i} \tau_{ik}G_i(\theta_i)}_{\sigma_{i4}}}
                + {\underbrace{%
                    \textstyle \sum_{k = 2}^{\hat{\rho}_i} \frac{\hat{a}_{ik}}{2\hat{b}_{ik}} \tau_{ik}^2}_{\sigma_{i5}}} ,
            \\
             \check{q}_j(\zeta_j ,\check{R}_j) & =
            {\underbrace{%
                \textstyle \sum_{k =2}^{\check{\rho}_j} \frac{\check{a}_{jk}}{2\check{b}_{jk}} \zeta_{jk}^2}_{\sigma_{j1}}},
        \end{aligned}
    \end{equation}
\end{figure*}
\end{proposition}}
\begin{remark}
\li{Proposition \ref{prop_main_result} shows that the algorithm dynamics generated via augmented Lagrangian methods can be given an optimal control interpretation. In particular, a cost functional is minimized penalizing the deviation of trajectories from the optimal equilibrium point during the algorithm transients. It should be noted that the minimization in \eqref{eq_primal_cost} is over all \lhi{stabilizing} \ils{control} inputs $U$, i.e., at any time~$t$, $U(t)$ is allowed to be an arbitrary functional of the history of the process up to that time, \lic{such that \licc{variables $(\theta,\lambda,\mu$)} 
converge to the optimal solution \lhj{$(\theta^\ast,\lambda^\ast,\mu^\ast)$, which is the same solution} the primal-dual dynamics \eqref{eq_comp_dynamics} converge to}. The significance of the cost functional used in \eqref{eq_primal_cost} and its interpretation will be discussed in Section \ref{sec_discussion}.}
\end{remark}

\subsection{Intuition}
\li{A detailed derivation of Proposition \ref{prop_main_result} is provided in Section \ref{sec_proof_p1}.
Here we provide some intuition for the derivation.}

\lhb{Typically, the solutions to optimal control problems of the form \eqref{eq_primal_cost} are 
\li{quantified via}
the Hamilton-Jacobi-\lhe{Bellman} equation \cite{liberzon_CalculusVariations_11,bertsekas_DynamicProgramming_00}, which is derived \ils{via dynamic programming arguments}. This \ils{is a} 
\lhk{\ilc{partial} differential equation} \ils{that is in general difficult to solve analytically} for nonlinear \ilc{systems.}}

\lhb{The unique solution to the Hamilton-Jacobi\lhk{-Bellman} equation \ils{is the \textit{value function} for the problem, which 
provides the optimal value of the cost considered,} but can also serve as a Lyapunov function.
Correspondingly, methods \ilc{exist, known as \textit{inverse optimal control} \cite{kalman_WhenLinear_64, moylan_NonlinearRegulator_73, sepulchre_ConstructiveNonlinear_97},} which link Lyapunov functions to system performance via an optimal control problem. The systems \lhl{analyzed} in \li{more classical inverse optimal control approaches}
\ils{involve centralized control implementations rather than}
a distributed network system, such as the one considered here.}


\lhb{\ilc{For 
network systems as in Figure \ref{fig_dist_opt} comprised of node and edge subsystems,
a network Lyapunov function can be constructed when these subsystems are passive, via
their local 
storage functions.}
In our case, for dynamics \eqref{eq_comp_dynamics} and \eqref{eq_comp_outputs}, the \lhn{\ilc{subsystems}} $\Sigma_{\theta i}$, $\Sigma_{\lambda i}^+$ and $\Sigma_{\mu j}$ \ilc{are passive with} storage functions $V_i(\xi_i): \mathbb{R}^{\rho_i} \rightarrow \mathbb{R}_{\geq 0}$, $\hat{V}_i(\tau_i) : \mathbb{R}^{\hat{\rho}_i} \rightarrow \mathbb{R}_{\geq 0}$ and $\check{V}_j(\zeta_j) : \mathbb{R}^{\check{\rho}_j} \rightarrow \mathbb{R}_{\geq 0}$ given by \cite{yamashita_PassivitybasedGeneralization_20}:
\begin{subequations} \label{eq_storage_local}
    \begin{align}
        V_i(\xi_i) & =  \textstyle \frac{1}{2b_{i1}} \left(\xi_{i1} - \theta_i^{\ast} \right)^2 + \sum_{k = 2}^{\rho_i} \frac{1}{2b_{ik}} \xi_{ik}^2, \label{eq_storage_local_1} \\
        \hat{V}_i(\tau_i) & = \textstyle \frac{1}{2 \hat{b}_{i1}} (\tau_{i1} - \lambda_i^{\ast})^2 + \sum_{k = 2}^{\hat{\rho}_i} \frac{1}{2 \hat{b}_{ik}} \tau_{ik}^2, \label{eq_storage_local_2} \\
        \check{V}_j(\zeta_j) & = \textstyle \frac{1}{2 \check{b}_{j1}} \left( \zeta_{j1} - \mu_j^{\ast} \right)^2 + \sum_{k = 2}^{\check{\rho}_j} \frac{1}{2 \check{b}_{jk}} \zeta_{jk}^2, \label{eq_storage_local_3}
    \end{align}
\end{subequations}
and the system Lyapunov function is given by
\begin{equation} \label{eq_storage_full}
    V(\xi,\tau,\zeta)  = \sum_{i \in \calv} V_i(\xi_i) + \sum_{i \in \calv} \hat{V}_i(\tau_i) + \sum_{j \in \cale} \check{V}_j(\zeta_j).
\end{equation}
}

We showed in our previous work \cite{hallinan_InverseOptimal_23a} that an inverse optimal control interpretation is possible for some \ilc{classes of} passive network systems, provided a set of sufficient conditions are satisfied.
\lhb{This result, combined with those of \cite{yamashita_PassivitybasedGeneralization_20}, indicates that \lic{an} optimal control interpretation of the algorithm \licc{in} \eqref{eq_comp_outputs} may be possible. However, neither the results of \cite{hallinan_InverseOptimal_23a}, nor standard results in inverse optimal control apply to systems like \eqref{eq_comp_dynamics} due to the non-smooth nature of the system $\Sigma_\lambda^+$, and the fact that it is non-affine in its input $\theta$, as shown in Figure \ref{fig_dist_opt}.}

\lhb{Despite this, we show \li{in this paper} 
through a first principles derivation that an optimal control interpretation is still possible for the augmented algorithm \eqref{eq_comp_dynamics} by using the sum of storage functions as a value function for this system. The proof, found in Section \ref{sec_proof_p1}, 
involves directly evaluating the cost functional \licc{in} \eqref{eq_primal_cost}
\lic{and exploiting the structure of the algorithm dynamics to derive various intermediate results that allow to 
show that the minimum is achieved when  $U = \underline{U}$.}}
\subsection{Discussion} \label{sec_discussion}
As 
\ilo{demonstrated  
in Section \ref{sec_improved_performance} via examples,} 
augmenting the dynamics 
as in
\eqref{eq_comp_dynamics}
can 
enhance transient performance.
\ilo{In this section we discuss that the cost functional \eqref{eq_primal_cost}, can be 
seen as penalizing transient deviations from optimality. It can therefore serve as a network performance metric that can be also} \lha{used to aid the selection of the algorithm parameters.}


\lhb{In order to interpret the cost function in \eqref{eq_primal_cost}, we use the following results,} \lho{which are proved in Appendix \ref{sec_proof_lemmas}}.



\begin{lemma} \label{lemma_state_cost_bound}
    \lhb{
    \lic{The terms defined in
    \eqref{eq_state_cost_primal_dual}, i.e., the terms in the  cost function \eqref{eq_primal_cost} involving the system states,
    satisfy the following relationships:} 
    \begin{subequations}
        \begin{align}
             \sigma_{i1} & = D_{F_i}(\theta_i,\ils{\theta_i^{*}}) + D_{F_i}(\ils{\theta_i^*},\theta_i), \label{eq_term_interp1} \\
             \begin{split}
                \sigma_{i2} + \sigma_{i4}  & \geq
                D_{G_i} (\theta_i^\ast,\theta_i)\lambda_i
                + D_{G_i} (\theta_i,\theta_i^{\ast}) \lambda_i^\ast
                \\ & \qquad -  G_i(\theta_i^\ast) \lambda_i,
                \label{eq_term_interp2}
            \end{split}
        \end{align}
    \end{subequations}
    where $D_{F_i}(\cdot,\cdot)$ and $D_{G_i}(\cdot,\cdot)$ \li{is the Bregman divergence} associated with $F_i(\cdot)$ and $G_i(\cdot)$ respectively, defined in \eqref{eq_bregman_def}.}
\end{lemma}

\begin{lemma} \label{lemma_non_neg}
    \lhb{The \lic{integrand in the cost functional} in \eqref{eq_primal_cost} is non-negative. Moreover, it \li{is equal} 
    to zero when the algorithm \licc{in} \eqref{eq_primal_dual_control} with optimal controller \eqref{eq_augmented_controllers} \li{is at} 
    equilibrium.}
\end{lemma}

\lhb{Lemma \ref{lemma_non_neg} shows that the optimal control problem is well-posed,
\lic{with the integrand in the cost functional tending to zero as the states tend to their equilibrium values.}
We observe that the costs associated with the controllers in \eqref{eq_primal_dual_R} and the costs associated with the corresponding states ($\sigma_{i3}$, $\sigma_{i5}$ and $\sigma_{j1}$ in \eqref{eq_state_cost_primal_dual}) are inversely related by the parameters $a_{ik}$, $\hat{a}_{ik}$ and $\check{a}_{jk}$, $k \geq 2$, which demonstrates a trade-off between control effort and state regulation, a property often seen in optimal control problems (e.g., in linear quadratic regulator (LQR) problems \cite{liberzon_CalculusVariations_11}).}

\lhb{Via \eqref{eq_term_interp1} and \eqref{eq_term_interp2}}, we can see that \lhk{$q_i(\xi_i,\tau_i,R_i)$ in} \eqref{eq_state_cost_primal_dual} introduces a cost associated with the deviation 
$\theta_i$ from the equilibrium point $\ils{\theta_i^{*}}$ via
\lhb{the Bregman divergences associated with $F_i(\cdot)$ and $G_i(\cdot)$. The Bregman divergence can be interpreted as a measure of distance using convex functions \cite{censor_ParallelOptimization_97,beck_MirrorDescent_03,banerjee_ClusteringBregman_05}, which demonstrates that the cost function in \eqref{eq_primal_cost} penalizes 
\lic{deviations in the primal variable from its equilibrium value,} while the algorithm is in the transient stage.}

\lhb{\ilo{Furthermore}, the terms $\sigma_{i3}$, $\sigma_{i5}$ and $\sigma_{j1}$ are associated with the additional auxiliary variables introduced by the augmented algorithm. The costs associated with these states are quadratic, and}
penalize the deviation of these variables from zero, which is their desired steady state value for the optimal solution of \eqref{opt_distributed_transformed} to be reached.

\lhb{\ilo{In addition,} we see in \eqref{eq_term_interp2} that $\sigma_{i2}+\sigma_{i4}$ is bounded below by $-G_i(\theta_i^\ast)\lambda_i$, which \ils{is a} positive quantity (as $G_i(\theta_i^\ast) \leq 0$ by the constraint in \eqref{opt_distributed} and $\lambda_i \geq 0$ by the operator \eqref{eq_positive_projection}) and is equal to zero when $\lambda_i = \lambda_i^\ast$ by the KKT condition \eqref{eq_distributed_kkt_2}). If $G_i(\theta_i^\ast) \neq 0$, then $\lambda_i^\ast =0$ and we can see that this bound linearly penalizes deviations of $\lambda_i$ from its equilibrium value.}

\lhb{Finally, we can observe \li{various} 
interesting features of the structure of the optimal control problem \eqref{eq_primal_cost}. Firstly, it can be seen that the cost function is separable into terms that correspond to each of the local network \li{subsystems.}
This means the optimal control problem can easily be scaled to account for new nodes in the network. Secondly,}
\lha{we also note that
the equilibrium point in the optimal control formulation \li{can be} 
arbitrary, and the result holds for any
\lhb{optimal solution to the optimization problem \eqref{opt_distributed_transformed}.}}

\vspace{-.05cm}
\section{Further Results} \label{sec_further_results}
\lhb{We now highlight some further refinements and extensions of the results presented in Section \ref{sec_opt_control_interp}. First, \li{we consider the special case}
where there are no constraints in \eqref{opt_distributed}. Next, we consider the case where the consensus constraint in \eqref{opt_distributed_transformed} is replaced by 
\lic{coupling inequality constraints}. And finally, we \ils{consider} 
the effect of including feed-forward terms in the system outputs in \eqref{eq_comp_outputs}.}

\subsection{Unconstrained Consensus Problem with No Feed-Forward Terms} \label{sec_unconstrained_no_ff}
\lhb{The results of Section \ref{sec_opt_control_interp} greatly simplify in the case where}
\lha{
there are no inequality constraints in \eqref{opt_distributed_transformed}.
\lhk{In this case, the dynamics associated with the dual variable $\lambda$ do not appear in the primal-dual algorithm \eqref{eq_primal_dual_alg} or the augmented algorithm \eqref{eq_comp_dynamics} (i.e., $\Sigma_{\lambda i}^+ = 0$ in \eqref{eq_comp_dynamics} and Figure \ref{fig_dist_opt}).}
\lhb{Again, we consider} the case where no feed-forward terms are used in \eqref{eq_comp_outputs} (i.e., $d_i=0$ and $\check{d}_j=0$ for all $i \in \calv$ and $j \in \calv$).}
\lhk{In control form, we can write the reduced dynamics as
\begin{subequations}
    \begin{align} \label{eq_primal_dual_control_uncon}
        &\Sigma_{\theta i} = \begin{cases}
            \dot{\xi}_i = b_i \mathbb{1}_{\rho_i} \phi_i + b_{ui} u_i, \\
            \theta_i = \mathbb{1}_{\rho_i}^T \xi_i,
        \end{cases} \\
        &\Sigma_{\mu j} = \begin{cases}
            \dot{\zeta}_j = \check{b}_j \mathbb{1}_{\check{\rho}_j} \omega_j + \check{b}_{uj} \check{u}_j, \\
            \mu_j = \mathbb{1}_{\check{\rho}_j}^T \zeta_j,
        \end{cases}
    \end{align}
\end{subequations}
where parameters and variables have equivalent definitions to those in \eqref{eq_primal_dual_control}.
Now, we recover the augmented algorithm \eqref{eq_comp_dynamics} if, for each $i \in \calv$ and $j \in \cale$, we let $u_i$ and $\check{u}_j$ equal to
\begin{subequations}
    \begin{align} \label{eq_augmented_controllers_uncon}
        \underline{u}_i &= \left[ -a_{i2} \xi_{i2}, \ldots, -a_{i\rho_{i}} \xi_{i\rho_{i}} \right]^T,    \\
        \underline{\check{u}}_j &= \left[ -\check{a}_{j2} \zeta_{j2}, \ldots, -\check{a}_{j\check{\rho}_j} \zeta_{j\check{\rho}_j} \right]^T.
    \end{align}
\end{subequations}
}%
\li{The} following result is \lhk{then} a direct corollary\footnote{\ilc{It should be noted that this corollary follows also from the results in \cite{hallinan_InverseOptimal_23a} for smooth systems that are affine in the input.}} of Proposition \ref{prop_main_result} when the dual dynamics associated with the inequality constraints are ignored, and the proof follows the same approach given in Section \ref{sec_proof_p1}, \lhk{with similar interpretation to that given in Section \ref{sec_discussion}}.

\begin{corollary} \label{cor_unconstrained_no_ff}
    Consider the following optimal control problem for \eqref{eq_primal_dual_control_uncon} with initial condition $(\xi_0, \zeta_0)$:
    \begin{equation} \label{eq_primal_cost_uncon}
        \min_{U} \int_0^{\infty} \lVert U \rVert_R^2 + \sum_{i \in \calv} q_i(\xi_i,R_i) + \sum_{j \in \cale} \check{q}_j(\zeta_j,\check{R}_j) \ dt,
    \end{equation}
    where $U = \begin{bmatrix} [u_i]_{i \in \calv}^T & [\check{u}_j]_{j \in \cale}^T \end{bmatrix}^T$. Then $\underline{U} = \begin{bmatrix} [\underline{u}_i]_{i \in \calv}^T & [\underline{\check{u}}_j]_{j \in \cale}^T \end{bmatrix}^T$ given by \eqref{eq_augmented_controllers_uncon} is the optimal solution to the optimal control problem \licc{in} \eqref{eq_primal_cost_uncon}, with \[R= \mathrm{diag} [\oplus_{i \in \calv} R_i ,\ \oplus_{j \in \cale} \check{R}_j ], \] where
    \begin{subequations}
        \begin{align}
            R_i &=  \textstyle \frac{1}{2} \mathrm{diag} \left[ \frac{1}{a_{i2}b_{i2}}, \ldots, \frac{1}{a_{i\rho_i}b_{i\rho_i}} \right], \\
            \check{R}_j &= \textstyle \frac{1}{2} \mathrm{diag} \left[ \frac{1}{\check{a}_{j2}\check{b}_{j2}}, \ldots, \frac{1}{\check{a}_{j\check{\rho}_j}\check{b}_{j\check{\rho}_j}} \right],
        \end{align}
    \end{subequations}
    and
    \begin{subequations} \label{eq_state_cost_primal_dual_uncon}
        \begin{align}
            &\begin{aligned}
             q_i(\xi_i,R_i)   &=  \left(\theta_i - \theta_i^{\ast} \right) \left(\nabla F_i(\theta_i) - \nabla F_i(\theta_i^{\ast}) \right)
                \\&\qquad + \textstyle \sum_{k = 2}^{\rho_i} \frac{a_{ik}}{2b_{ik}} \xi_{ik}^2,
            \end{aligned} \label{eq_state_cost_primal_uncon} \\
            &\check{q}_j(\zeta_j ,\check{R}_j)   =  \textstyle \sum_{k =2}^{\check{\rho}_j} \frac{\check{a}_{jk}}{2\check{b}_{jk}} \zeta_{jk}^2. \label{eq_state_cost_dual_uncon}
        \end{align}
    \end{subequations}
\end{corollary}



\subsection{
\lic{\lhi{Coupling} Inequality Constraints}} \label{sec_global_ineq}
\lhb{In some applications, e.g., communication networks \cite{kelly_RateControl_98}, the optimization problem does not take the form of a consensus problem like \eqref{opt_distributed}, but instead uses inequality constraints that act on a set of primal variables associated with the agents. }

\lhb{Such problems can be \lic{formulated} 
as follows:
\begin{mini}
	{\theta \in \mathbb{R}^{|\calv|}}{\sum_{i \in \calv} F_i(\theta_i),}{\label{opt_dist_ineqality}}{}
        \addConstraint{H(\omega) & \leq 0, }
\end{mini}
where \lhi{here} $\omega := [\omega_j]_{j\in\cale} = \lhi{\mathcal{R}} \theta \in \mathbb{R}^{|\cale|}$,
\lhi{$\mathcal{R}\in \mathbb{R}^{|\cale|\times |\calv|}$
\lic{and} }
$H_j: \mathbb{R} \rightarrow \mathbb{R}$ is convex and continuously differentiable.}
\lhk{Again, we assume Slater's condition holds and  $\sum_{i\in\calv} F_i(\theta_i)$ is bounded from below in the feasible set of \eqref{opt_dist_ineqality}.}
\begin{remark}
    \lhi{
    \lic{The notation above, where the argument of $H(.)$ is a} linear \licc{combination} of the primal variables via $\omega = \mathcal{R}\theta$, is \lic{used} 
    as it appears \lic{often} 
    in applications. For example, $\mathcal{R}$ can be a routing matrix in network congestion control and multi-path routing problems, and the constraints $H(\omega)\leq 0$ ensure that the sum of data transmission rates along a link in the network is less than or equal to the link's capacity \cite{kelly_RateControl_98, holding_StabilityInstability_15}.}
\end{remark}

\begin{figure}[t]
\centering
\small
\input{Images/tikz_inequality}
\caption{\lhb{Closed loop block diagram of the primal-dual distributed optimization algorithm \licc{in} \eqref{eq_primal_dual_ineq} used to solve a problem of the form \eqref{opt_dist_ineqality} with a distributed inequality constraint. The primal dynamics associated with the nodes are given by system $\Sigma_\theta$ and dual dynamics for the inequality constraints associated with the edges are given by $\Sigma_\lambda^+$}.}
\label{fig_dist_ineq}
\end{figure}
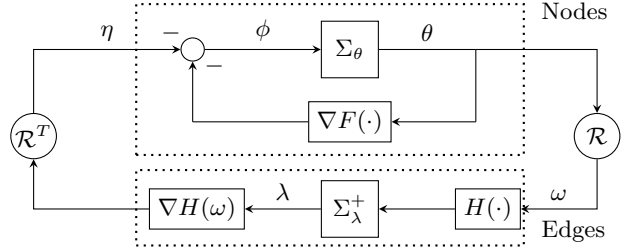

The KKT conditions associated with \eqref{opt_dist_ineqality} are
\begin{subequations} \label{eq_ineq_kkt}
    \begin{align}
        & \nabla F(\theta^{\ast}) + \lhi{\mathcal{R}^T} \nabla H(\omega) \lambda^{\ast}  = 0 \label{eq_ineq_kkt_1} \\
        & \lambda_j^{\ast}H_j(\omega_j^{\ast}) = 0, \quad \lambda_j^{\ast} \geq 0, \quad H_j(\omega_j^{\ast}) \leq 0,  \label{eq_ineq_kkt_2}
    \end{align}
\end{subequations}
where $\lambda^{\ast} := [\lambda_j^{\ast}]_{i \in \cale} \in \mathbb{R}^{|\cale|}$ is an optimal value of the dual variable $\lambda := [\lambda_j]_{j \in \cale} \in \mathbb{R}^{|\cale|}$ associated with the inequality constraints, $F(\theta) = \sum_{j \in \cale} F_i(\theta_i)$, and \lhj{$\nabla H(\omega) = \mathrm{diag}[ \nabla H_j(\omega_j)]_{j\in\cale}$}, and \eqref{eq_ineq_kkt_2} holds for all $ j \in \cale$.
\lhb{\lhk{These admit the following primal-dual dynamics,  where the optimal solution $(\theta^\ast,\lambda^\ast)$ is reached at equilibrium:} 
\begin{subequations} \label{eq_primal_dual_ineq}
    \begin{align}
        \dot{\theta} &= -\nabla F(\theta) - \eta, \\
        \dot{\lambda} & = [H(\omega)]_{\lambda}^+, \ \lambda(0) \geq 0,
    \end{align}
\end{subequations}
where now
\begin{equation} \label{eq_def_eta_ineq}
    \eta= [\eta_i]_{i \in \calv} := \lhi{\mathcal{R}^T} \nabla H(\omega) \lambda,
\end{equation}
and, again, we use the operator \eqref{eq_positive_projection}.
The dynamics are summarized in the block diagram in Figure \ref{fig_dist_ineq}, with primal dynamics represented by $\Sigma_\theta$, dual dynamics represented by $\Sigma_\lambda^+$ and $\phi = [\phi_i]_{i \in \calv} := - \nabla F(\theta) - \eta$.
We see that when compared to Figure \ref{fig_dist_opt}, the dual dynamics associated with the inequality constraints are represented by the edges of the network.}

Now, we can augment the dynamics \licc{in} \eqref{eq_primal_dual_ineq}, where each $\Sigma_{\theta i}$ and $\Sigma_{\lambda j}^+$ in Figure \ref{fig_dist_ineq} is replaced by:
\begin{subequations}\label{eq_comp_dynamics_ineq}
    \begin{align}
        & \Sigma_{\theta i} : \begin{cases}
            \dot{\xi}_{i1} &= b_{i1} \phi_i, \\
            \dot{\xi}_{ik} &= b_{ik} \phi_i - a_{ik} \xi_{ik}, \quad k = 2, \ldots, \rho_i,
        \end{cases} \label{eq_comp_dynamics_xi_ineq} \\
        & \Sigma_{\lambda j}^+ : \begin{cases}
            \dot{\tau}_{j1} &= [\hat{b}_{j1} H_j(\omega_j)]_{\tau_{j1}}^+, \\
            \dot{\tau}_{jk} &= [\hat{b}_{jk} H_j(\omega_j) - \hat{a}_{jk} \tau_{jk}]_{\tau_{jk}}^+, \quad k = 2, \ldots, \hat{\rho}_j,
        \end{cases} \label{eq_comp_dynamics_tau_ineq}
    \end{align}
\end{subequations}
and the outputs of $\Sigma_{\theta i}$ and $\Sigma_{\lambda j}^+$ are
\begin{align}
    &\theta_i = \textstyle \sum_{k = 1}^{\rho_i} \xi_{ik},
    & \lambda_j = \textstyle \sum_{k = 1}^{\hat{\rho}_j} \tau_{jk}.
\end{align}
Here, the parameters are defined similarly to those in \eqref{eq_comp_dynamics}. We can reformulate these dynamics in control form as follows:
\begin{subequations}
    \begin{align}  \label{eq_primal_dual_control_ineq}
        &\Sigma_{\theta i} = \begin{cases}
            \dot{\xi}_i = b_i \mathbb{1}_{\rho_i} \phi_i + b_{ui} u_i, \\
            \theta_i = \mathbb{1}_{\rho_i}^T \xi_i,
        \end{cases} \\
        &\Sigma_{\lambda j}^+ = \begin{cases}
            \dot{\tau}_j = [ \hat{b}_j \mathbb{1}_{\hat{\rho}_j} H_j(\omega_j) + \hat{b}_{uj} \hat{u}_j ]^+_{\tau_j}, \\
            \lambda_j = \mathbb{1}_{\hat{\rho}_j}^T \tau_j,
        \end{cases}
    \end{align}
\end{subequations}
where the parameters follow similar definition to those in \eqref{eq_primal_dual_control}. We see that we recover the augmented algorithm \eqref{eq_comp_dynamics_ineq} if, for each $i \in \calv$ and $j \in \cale$, we set each $u_i$ and $\hat{u}_j$ in \eqref{eq_primal_dual_control_ineq} equal to
\begin{subequations}
    \begin{align}  \label{eq_augmented_controllers_ineq}
        \underline{u}_i &= \left[ -a_{i2} \xi_{i2}, \ldots, -a_{i\rho_{i}} \xi_{i\rho_{i}} \right]^T,  \\
        \underline{\hat{u}}_j &= \left[ -\hat{a}_{j2} \tau_{j2}, \ldots, -\hat{a}_{j\hat{\rho}_{j}} \tau_{j\hat{\rho}_{j}} \right]^T.
    \end{align}
\end{subequations}

We now see in Corollary \ref{cor_inequality} that the dynamics \eqref{eq_comp_dynamics_ineq} can be obtained via an optimal control formulation.

\begin{corollary} \label{cor_inequality}
    Consider the following optimal control problem for \eqref{eq_primal_dual_control_ineq} with initial condition $(\xi_0, \tau_0)$:
    \begin{equation*} 
        \min_{U} \int_0^{\infty} \lVert U \rVert_R^2 + \sum_{i \in \calv} q_i(\xi_i,\tau,R_i) + \sum_{j \in \cale} \hat{q}_j(\tau_j,\xi,\check{R}_j) \ dt,
    \end{equation*}
    where $U = \begin{bmatrix} [u_i]_{i \in \calv}^T & [\hat{u}_j]_{j \in \cale}^T \end{bmatrix}^T$. Then $\underline{U} = \begin{bmatrix} [\underline{u}_i]_{i \in \calv}^T & [\underline{\hat{u}}_j]_{j \in \cale}^T \end{bmatrix}^T$ given by \eqref{eq_augmented_controllers_ineq} is the optimal solution to the optimal control problem
    , with $R= \mathrm{diag} [\oplus_{i \in \calv} R_i , \oplus_{j \in \cale} \hat{R}_j ]$, where
    \begin{subequations}
        \begin{align} \label{eq_primal_dual_R_ineq}
            R_i &= \textstyle \frac{1}{2} \mathrm{diag} \left[ \frac{1}{a_{i2}b_{i2}}, \ldots, \frac{1}{a_{i\rho_i}b_{i\rho_i}} \right], \\
            \hat{R}_j &= \textstyle \frac{1}{2} \mathrm{diag} \left[ \frac{1}{\hat{a}_{j2}\hat{b}_{j2}}, \ldots, \frac{1}{\hat{a}_{j\hat{\rho}_j}\hat{b}_{j\hat{\rho}_j}} \right],
        \end{align}
    \end{subequations}
    and
     \begin{align*}
            & \begin{aligned}
            q_i(\xi_i,\tau,R_i)   &=  \left(\theta_i - \theta_i^{\ast} \right) \left(\nabla F_i(\theta_i) - \nabla F_i(\theta_i^{\ast})
                + \eta_i - \eta_i^{\ast} \right)
                \\ & \qquad + \textstyle \sum_{k = 2}^{\rho_i} \frac{a_{ik}}{2b_{ik}} \xi_{ik}^2,
            \end{aligned} \\
            & \begin{aligned}
            \hat{q}_j(\tau_j,\xi,\hat{R}_j)  & = -(\tau_{j1} - \lambda_j^{\ast})[H_j(\omega_j)]_{\tau_{j1}}^+
            \\& \qquad  - \textstyle \sum_{k = 2}^{\hat{\rho}_j} \tau_{jk}H_j(\omega_j)
            \textstyle+ \sum_{k = 2}^{\hat{\rho}_j} \frac{\hat{a}_{jk}}{2\hat{b}_{jk}} \tau_{jk}^2.
            \end{aligned}
        \end{align*}
\end{corollary}

\lhb{
The proof directly follows the approach taken in Section \ref{sec_proof_p1} for Proposition 1, but we note here that a similar relation to \eqref{eq_proof_relation} cannot be derived for this system, \li{and} 
Step 1 of the proof can be omitted.
Unlike the state cost function
\eqref{eq_state_cost_primal_dual} for the case \lhl{analyzed} in Proposition \ref{prop_main_result}, we see that the state cost functions now depends on the network interconnection matrix $\cala$ via the definition of $\eta_i$ in \eqref{eq_def_eta_ineq}. However, the structure of the state cost function is otherwise \lic{analogous} 
to that in \eqref{eq_state_cost_primal_dual}, and the terms have the same interpretation. }

\subsection{Unconstrained Consensus Problem with Feed-Forward} \label{sec_feedforward}
Here, we present a modification of the results 
\lhk{of Section \ref{sec_opt_control_interp}} to account for input feed-forward terms at the edges. \lha{\lhb{We will see that this will allow} us to determine a performance metric for algorithms of the form \eqref{eq_primal_dual_lagrangian}}.
\lha{\lhk{As in Section \ref{sec_unconstrained_no_ff}}, 
we consider the unconstrained case of \eqref{opt_distributed_transformed} (i.e., $\Sigma_\lambda^+ = 0$ in Figure \ref{fig_dist_opt}) in order to simplify the presentation, but an analogous approach to Proposition \ref{prop_main_result} can be followed to derive a performance result for optimization problems involving inequality constraints.}

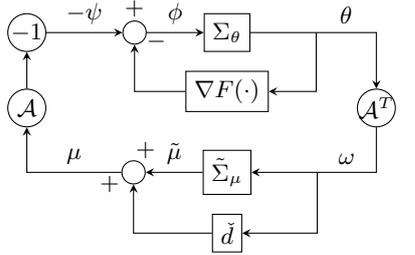
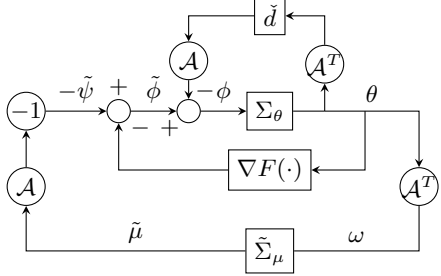
\begin{figure}[t]
\centering
\begin{subfigure}[b]{0.45\textwidth}
    \centering
    \input{Images/tikz_feed_forward}
    \caption{Input feed-forward terms at edges.}
    \label{fig_feedforward}
\end{subfigure}
%
%
\begin{subfigure}[b]{0.45\textwidth}
    \centering
    \input{Images/tikz_feedback}
    \caption{Output feedback terms at nodes. }
    \label{fig_feedback}
\end{subfigure}
\caption{Reformulation of Figure \ref{fig_dist_opt} demonstrating input feed-forward terms at the edges transformed to output feedback at the nodes \lhf{for the unconstrained problem (i.e., $\Sigma_\lambda^+ = 0$)}.}
\label{fig_transformation}
\end{figure}

Consider the \lhm{subsystem outputs \eqref{eq_comp_outputs_x} and \eqref{eq_comp_outputs_mu}} with $d_i = 0, \forall i \in \calv$ and $\check{d}_j \geq 0, \forall j \in \cale$. The extra term $\check{d}_j\omega$ at each edge $j \in \cale$ acts as an input feed-forward term,
\lhk{as} illustrated in Figure \ref{fig_feedforward}. Here, we define $\check{d} = \oplus_{j \in \cale} \check{d}_j$, $\tilde{\mu}_j =\mathbb{1}_{\check{\rho}_j}^T \zeta_j$, $\tilde{\mu} = [\tilde{\mu}_j]_{j \in \cale}$, and \lhm{$\tilde{\Sigma}_\mu$ is the block-diagonal edge system with input $\omega$ and output $\tilde{\mu}$ (i.e., the system $\Sigma_\mu$} without the input feed-forward term, see Figure \ref{fig_feedforward}). \lha{We also take $\cala$ to be the directed network incidence matrix.}

\lha{In order} to apply
\lha{the results of \cite{hallinan_InverseOptimal_23a} to this system,}
the input feed-forward term at the edges must be moved to the nodes, which is demonstrated in Figure \ref{fig_feedback}. Here, we denote
\begin{subequations} \label{eq_def_psi_phi_tilde}
    \begin{align}
        \tilde{\psi} &= \cala \tilde{\mu}, &
        \tilde{\phi} &= -\nabla F(x) - \tilde{\psi},
    \end{align}
\end{subequations}
\lhe{and all other variables are as in \eqref{eq_var_defs}.}

With this transformation, the composite node and edge subsystems in Figure \ref{fig_feedback} can be written as
\begin{subequations} \label{eq_opt_ext_dynamics}
    \begin{align}
        & \dot{\xi}  = b \mathbb{1}_\calv \tilde{\phi} - \left( I_{o\calv} a I_{o\calv}^T  + b \mathbb{1}_\calv \mathcal{L}_{\check{d}} \mathbb{1}_\calv^T \right) \xi \\
        & \dot{\zeta} = \check{b} \mathbb{1}_{\cale} \lhe{\omega} - I_{o \cale} \check{a} I_{o \cale}^T \zeta
    \end{align}
\end{subequations}
where we define $\xi = [\xi_i]_{i \in \calv}$, $b = \oplus_{i \in \calv} b_i$, $ a = \oplus_{i \in \calv} a_i$ \lhb{(with $a_i = \mathrm{diag} [a_{i2}, \ldots, a_{i\rho_i}]$)}, $I_{o\calv} = \oplus_{i \in \calv} I_{o \rho_i}$ and $\mathbb{1}_\calv = \oplus_{i \in \calv} \mathbb{1}_{\rho_i}$, with similar definitions holding for $\zeta$, $\check{b}$, $\check{a}$, $I_{o \cale}$ and $\mathbb{1}_\cale$.
\lha{Here, $\mathcal{L}_{\check{d}} := \cala \check{d} \cala^T$ is the network Laplacian weighted by matrix $\check{d}$.}

We now follow a similar approach to that taken in Section \ref{sec_unconstrained_no_ff}. The dynamics \eqref{eq_opt_ext_dynamics} can be written in control form as
\begin{subequations}
    \begin{align} \label{eq_primal_dual_control_uncon_ff}
        &\Sigma_{\theta} = \begin{cases}
            \dot{\xi} = b \mathbb{1}_\calv \tilde{\phi} + b_{u\calv} u_\calv,  \\
            \theta = \mathbb{1}_{\calv}^T \xi,
        \end{cases} \\
        &\tilde{\Sigma}_{\mu} = \begin{cases}
            \dot{\zeta} = \check{b} \mathbb{1}_\cale \omega + b_{u\cale} u_\cale, \\
            \tilde{\mu} = \mathbb{1}_{\cale}^T \zeta,
        \end{cases}
    \end{align}
\end{subequations}
where \lhg{$\omega$ and $\tilde{\phi}$ are defined in \eqref{eq_def_z} and \eqref{eq_def_psi_phi_tilde} respectively, $\ils{u_{\calv}(t)} \in \mathbb{R}^{ \left( \Sigma_{i \in \calv} \left(\rho_i-1\right)\right) + |\cale|}$ and $\ils{u_{\cale}(t)} \in \mathbb{R}^{\Sigma_{j \in \cale} \left(\check{\rho}_j-1 \right)}$ are \lic{stabilizing} \ils{control} \lha{inputs}, and} we take
\begin{align}
    b_{u \calv} & = b \begin{bmatrix} I_{o\calv} & \mathbb{1}_\calv \cala \end{bmatrix}, &
    b_{u \cale} & = \check{b} I_{o\cale}.
\end{align}
We \lha{observe that we} recover the dynamics \eqref{eq_opt_ext_dynamics} if we use the \lha{distributed} controllers
\begin{align} \label{eq_opt_ext_controller}
    \underline{u}_\calv & = - \begin{bmatrix} a I_{o \calv}^T b^{-1} \\ \check{d} \cala^T \mathbb{1}_\calv^T \end{bmatrix} \xi, &
    \underline{u}_\cale & = - \check{a} I_{o \cale}^T \check{b}^{-1} \zeta.
\end{align}

\lhb{We now \lic{show} 
that the dynamics \eqref{eq_opt_ext_dynamics} can be constructed via an optimal control problem, providing a theoretical performance guarantee for such systems. The proof can be found in Appendix \ref{sec_proof_ff}.} 


\begin{proposition} \label{prop_feed_forward}
    \lha{Consider the following optimal control problem for system \eqref{eq_primal_dual_control_uncon_ff} with initial condition $(\xi_0, \zeta_0)$:}
    \begin{equation}  \label{eq_opt_ext_cost}
        \begin{aligned}
            \min_{U} \int_0^{\infty}  \lVert U \rVert_R^2  + \lhe{q_\calv(\xi,R_\calv) + q_\cale(\zeta,R_\cale)} \ dt
        \end{aligned}
    \end{equation}
    where $U = \begin{bmatrix} u_\calv & u_\cale \end{bmatrix}^T$.
    \lha{Then $\underline{U} = \begin{bmatrix} [\underline{u}_\calv & \underline{u}_\cale \end{bmatrix}^T$ given by \eqref{eq_opt_ext_controller} is the optimal solution to the optimal control problem \eqref{eq_primal_cost}, with}
    $R= \mathrm{blockdiag} [ R_\calv, R_\cale]$, where
    \begin{align} \label{eq_opt_ext_R}
        R_\calv^{-1} & = \begin{bmatrix} 2a  I_{o\calv}^T b^{-1} I_{o \calv} & 0 \\ 0 & 2 \check{d} \end{bmatrix}, & \!\!\!
        R_\cale^{-1} & =  2 \check{a}  I_{o\cale}^T \check{b}^{-1} I_{o \cale}.
    \end{align}
    and \lhe{
    \begin{subequations} \label{eq_state_cost_ff}
        \begin{align}
            q_\calv(\xi,R_\calv) & = \tfrac{1}{2}\theta^T \mathcal{L}_{\check{d}} \theta + \textstyle  \sum_{i \in \calv} q_i(\xi_i,R_i), \\
            q_\cale(\zeta,R_\cale) & = \textstyle  \sum_{j \in \cale} \check{q}_j(\zeta_j,\check{R}_j),
        \end{align}
    \end{subequations}
    with} $q_i(\xi_i,R_i)$, $\check{q}_j(\zeta_j ,\check{R}_j)$ given by \eqref{eq_state_cost_primal_dual_uncon}.
\end{proposition}

Here, we see that the cost functional \eqref{eq_opt_ext_cost} takes the same form as
\eqref{eq_primal_cost_uncon} in Corollary \ref{cor_unconstrained_no_ff}
\lha{(and as \eqref{eq_primal_cost} in Proposition \ref{prop_main_result}, without the terms corresponding to the dual dynamics for the inequality constraints)}, with an additional term \lha{
\begin{equation}
    \frac{1}{2}\theta^T \mathcal{L}_{\check{d}} \theta = \sum_{j \in \cale} \frac{1}{2} \check{d}_j \left(\theta_p - \theta_q \right)^2,
\end{equation}
where $p \in \calv$ and $q \in \calv$ denote the nodes connected by edge $j \in \cale$.
This term penalizes} the deviation from consensus in $\theta_i$ \lha{during algorithm transients, and is zero at equilibrium}.

\section{Conclusion} \label{sec_conclusion}
In this paper, we \icc{have}
\lha{
\li{provided} an optimal control interpretation \li{for}}
a general class of \li{augmented primal-dual distributed optimization algorithms.
\li{
In particular, we have shown 
that algorithms arising from augmented Lagrangian approaches 
can be seen as
minimizing a network-wide cost functional, which penalizes deviations from the optimal solution during the algorithm transients. The result has been shown for broad classes of dynamics, including the more involved setting where inequalities are present, and algorithms with feed-forward terms
that can further improve performance. The results derived in the paper provide 
an improved understanding of the performance enhancements augmented Lagrangian approaches provide, and more broadly improve our understanding of \lic{the} performance of distributed optimization algorithms.}}


\appendix
\section{Proof of Proposition \ref{prop_main_result}} \label{sec_proof_p1}

\lho{
\ilo{In this section we present the proof of} Proposition \ref{prop_main_result}.}
\lhh{Section \ref{sec_preliminaries} \lhj{first} provides some
preliminary definitions \li{of quantities used in the derivation of} 
our \ils{result.}
Next, Section \ref{sec_proof_p1_steps} provides the proof,}
\li{which is based \lic{on} inverse optimal control, i.e., approaches that link Lyapunov functions used to prove convergence, with cost functionals associated with performance. Classical approaches in inverse optimal control  
 are based on smooth dynamical systems,  which are affine in the input, and with also centralized control policies \cite{sepulchre_ConstructiveNonlinear_97}. The problem under consideration is more involved
as these features do not hold. In the derivation we show that the structure of the algorithm dynamics can be exploited to provide an optimal control interpretation despite the decentralized and nonlinear/non-smooth characteristics of the algorithm.}

\li{For the derivation we evaluate}
\lhh{the cost functional \licc{in} \eqref{eq_primal_cost}} under  \li{a \lhj{stabilizing}
} control action given by
\begin{equation} \label{eq_define_v}
    \begin{split}
        U  = \underline{U} + v,
    \end{split}
\end{equation}
where \lhe{$\underline{U} = \begin{bmatrix} [\underline{u}_i]_{i \in \calv}^T & \underline{\hat{u}}_i]_{i \in \calv}^T & [\underline{\check{u}}_j]_{j \in \cale}^T \end{bmatrix}^T$ is our proposed optimal controller 
\ilc{given by} \eqref{eq_augmented_controllers} and} $v$ is an \ils{arbitrary function of time} 
\ilc{
s.t. the system states converge to an equilibrium point where variables $(\theta,\lambda,\mu$) 
are equal to the optimal solution $(\theta^\ast,\lambda^\ast,\mu^\ast)$.}
\lhe{\li{It is then sufficient to show that} \eqref{eq_primal_cost} takes its minimum value when $v=0$, 
\li{as this then implies}
$\underline{U}$ is indeed the optimal controller.}

\subsection{Preliminary Definitions} \label{sec_preliminaries}

\lha{We first introduce some notation and definitions that will be used in the proof of Proposition \ref{prop_main_result}. Firstly, we let
\begin{align}
    R & = \mathrm{diag} \begin{bmatrix}
        \oplus_{i \in \calv} R_i,
        & \oplus_{i \in \calv} \hat{R}_i,
        & \oplus_{j \in \cale} \check{R}_j
    \end{bmatrix}, \label{eq_control_cost_ineq}
\end{align}
where
$R_i$, $\hat{R}_i$ and $\check{R}_j$ are given in \eqref{eq_primal_dual_R}. }
\lhe{Next, we define the cost functional to be \lhj{minimized} in \eqref{eq_primal_cost} as
    \begin{equation} \label{eq_J}
    \begin{split}
         J(U): = \int_0^{\infty} \lVert U \rVert_R^2 + q(\xi,\tau,\zeta,R) \ dt,
    \end{split}
\end{equation}
where $U = \begin{bmatrix} [u_i]_{i \in \calv}^T & [\hat{u}_i]_{i \in \calv}^T & [\check{u}_j]_{j \in \cale}^T \end{bmatrix}^T$, with $u_i$, $\hat{u}_i$ and $\check{u}_j$ control inputs for dynamics \eqref{eq_primal_dual_control}, $R$ is defined in \eqref{eq_control_cost_ineq}, and
\begin{equation} \label{eq_full_state_cost}
    \begin{split}
        q(\xi,\tau,\zeta,R) &= \sum_{i \in \calv} q_i(\xi_i,\tau_i, R_i)  + \sum_{i \in \calv} \hat{q}_i(\tau_i,\hat{R}_i) \\ &\qquad + \sum_{j \in \cale} \check{q}_j(\zeta_j,\check{R}_j)
    \end{split}
\end{equation}
with $q_i(\xi_i,\tau_i,R_i)$, $\hat{q}_i(\tau_i,\hat{R}_i)$ and $\check{q}_j(\zeta_j,\check{R}_j)$ defined in \eqref{eq_state_cost_primal_dual}.}
We now provide a proof of Proposition \ref{prop_main_result}.


\subsection{Proof} \label{sec_proof_p1_steps}
\lhb{We break the proof into \lhe{five} steps. In Step 1, we obtain an identity that relates the inputs and outputs of the node and edge subsystems. \lhe{In Step 2, we present some identities that relate the terms in our proposed optimal controller \licc{in} \eqref{eq_augmented_controllers} to the control cost matrix \eqref{eq_control_cost_ineq} and the storage functions \eqref{eq_storage_local}.} In Step \lhe{3}, we use \lhe{these identities} and the KKT conditions to re-write the state cost function \eqref{eq_full_state_cost} in terms of the local storage functions \eqref{eq_storage_local}, and in terms of our proposed optimal controller \eqref{eq_augmented_controllers}. \lhe{In Step 4, we derive an expression for the control cost $\lVert U \rVert_R$}. Finally, in Step \lhe{5}, we directly evaluate the cost functional \eqref{eq_J} using $U$ in \eqref{eq_define_v}, and show that the minimum value is reached when $v=0$, therefore \li{showing} 
that the optimal controller is $\underline{U}$.}
\\ \hbox{} \\
\textbf{Step 1:}
Firstly, we note that \lhe{
\begin{equation} \label{eq_storage_gradient}
    \begin{aligned}
    &\nabla V_i(\xi_i) = \begin{bmatrix}
        \frac{1}{b_{i1}} (\xi_{i1}-\theta_i^\ast), \frac{1}{b_{i2}} \xi_{i2}, \ \dots, \ \frac{1}{b_{i \rho_i}} \xi_{i \rho_i}
    \end{bmatrix}^T, \\
    &\nabla \hat{V}_i(\tau_i) = \begin{bmatrix}
        \frac{1}{\hat{b}_{i1}} (\tau_{i1}-\lambda_i^\ast), \frac{1}{\hat{b}_{i2}} \tau_{i2}, \ \dots, \frac{1}{\hat{b}_{i \hat{\rho}_i}} \tau_{i \hat{\rho}_i}
    \end{bmatrix}^T \\
    &\nabla \check{V}_i(\tau_i) = \begin{bmatrix}
        \frac{1}{\check{b}_{j1}} (\zeta_{j1}-\mu_j^\ast), \frac{1}{\check{b}_{j2}} \zeta_{j2}, \ \dots, \ \frac{1}{\check{b}_{j \check{\rho}_j}} \zeta_{j \check{\rho}_j}
    \end{bmatrix}^T, &
    \end{aligned}
\end{equation}
\lhk{where $V_i(\xi_i)$, $\hat{V}_i(\tau_i)$, and $\check{V}_i(\tau_i)$ are the subsystem storage functions defined in \eqref{eq_storage_local}}, from which we get
\begin{subequations} \label{eq_grad_identities}
    \begin{align}
        \begin{split}
            \nabla^T V_i(\xi_i) b_i \mathbb{1}_{\rho_i} & = \textstyle \sum_{k=1}^{\rho_i} \xi_{ik} - \theta_i^\ast  = \theta_i - \theta_i^{\ast},
        \end{split} \label{eq_grad_identities1} \\
        \begin{split}
            \nabla^T \hat{V}_i(\tau_i) \hat{b}_i \mathbb{1}_{\hat{\rho}_i} & = \textstyle \sum_{k=1}^{\hat{\rho}_i} \tau_{ik} - \lambda_i^\ast  = \lambda_i - \lambda_i^{\ast},
        \end{split} \label{eq_grad_identities2} \\
        \begin{split}
            \nabla^T \check{V}_j(\zeta_j) \check{b}_j \mathbb{1}_{\check{\rho}_j} & = \textstyle \sum_{k=1}^{\check{\rho}_j} \zeta_{jk} - \mu_j^\ast= \mu_j - \mu_j^{\ast}.
        \end{split} \label{eq_grad_identities3}
    \end{align}
\end{subequations}
If we aggregate all node and edge vectors together, we get}
\begin{subequations} \label{eq_comp_identity}
    \begin{align}
        [\nabla^T V_i(\xi_i) b_i \mathbb{1}_{\rho_i}]_{i \in \calv} & = \theta-\theta^{\ast}, \label{eq_comp_identity_1} \\
        [\nabla^T \hat{V}_i(\tau_i) \hat{b}_i \mathbb{1}_{\hat{\rho}_i}]_{i \in \calv} & = \lambda-\lambda^{\ast}, \label{eq_comp_identity_2} \\
        [\nabla^T \check{V}_j(\zeta_j) \check{b}_j \mathbb{1}_{\check{\rho}_j}]_{j \in \cale} & = \mu-\mu^{\ast}.
        \label{eq_comp_identity_3}
    \end{align}
\end{subequations}
By left-multiplying \eqref{eq_comp_identity_3} by $\cala$, and further left-multiplying by the transpose of \eqref{eq_comp_identity_1}, we arrive at the following \lha{relation using \eqref{eq_def_z} and \eqref{eq_def_psi}}:
\begin{align*}
    &([\nabla^T V_i(\xi_i) b_i\mathbb{1}_{\rho_i}]_{i \in \calv} )^T \cala ( \mu-\mu^{\ast} )
        \\ &\qquad = (\theta-\theta^{\ast})^T \cala ([\nabla^T \check{V}_j(\zeta_j) \check{b}_j \mathbb{1}_{\check{\rho}_j}]_{j \in \cale} ) \\
     \implies & ([\nabla^T V_i(\xi_i) b_i \mathbb{1}_{\rho_i}]_{i \in \calv})^T (\psi - \psi^{\ast})
     \\ & \qquad = (\omega - \omega^{\ast} )^T( [\nabla^T \check{V}_j(\zeta_j) \check{b}_j \mathbb{1}_{\check{\rho}_j}]_{j \in \cale} ).
\end{align*}
We therefore \ic{obtain} 
the following identity
\begin{equation} \label{eq_proof_relation}
    \begin{split}
        & \textstyle \sum_{i \in \calv} \nabla^T V_i(\xi_i) b_i \mathbb{1}_{\rho_i} (\psi_i - \psi_i^{\ast}) \\
         & \qquad- \textstyle \sum_{j \in \cale} \nabla^T \check{V}_j(\zeta_j) \check{b}_j \mathbb{1}_{\check{\rho}_j} (\omega_j - \omega_j^{\ast}) = 0.
     \end{split}
\end{equation}
\textbf{Step 2:}
\lhe{We now obtain some identities related to the terms in our proposed optimal control controller \eqref{eq_augmented_controllers}. Firstly we note that using \eqref{eq_storage_gradient}, we can write
\begin{subequations} \label{eq_identity_u0}
    \begin{align}
        & \underline{u}_i = \textstyle -\frac{1}{2}R_i^{-1}b_{ui}^T\nabla V_i(\xi_i), \\
        & \underline{\hat{u}}_i = \textstyle-\frac{1}{2}\hat{R}_i^{-1}\hat{b}_{ui}^T\nabla \hat{V}_i(\tau_i), \\
        & \underline{\check{u}}_j = \textstyle -\frac{1}{2}\check{R}_j^{-1}\check{b}_{uj}^T\nabla \check{V}_j(\zeta_j).
    \end{align}
\end{subequations}
Next, we see
\begin{subequations} \label{eq_identity_u1}
    \begin{align}
        & \lVert \underline{u}_i \rVert_{R_i}^2  = \textstyle \sum_{k =2}^{\rho_i} \frac{a_{ik}}{2b_{ik}} \xi_{ik}^2, \\
        & \lVert \underline{\hat{u}}_i \rVert_{\hat{R}_i}^2  = \textstyle  \sum_{k =2}^{\hat{\rho}_i} \frac{\hat{a}_{ik}}{2\hat{b}_{ik}} \tau_{ik}^2, \\
        & \lVert \underline{\check{u}}_j \rVert_{\check{R}_j}^2  = \textstyle \sum_{k =2}^{\check{\rho}_j} \frac{\check{a}_{jk}}{2\check{b}_{jk}} \zeta_{jk}^2.
    \end{align}
\end{subequations}
}%
\textbf{Step 3:}
\lhe{Using \eqref{eq_grad_identities} and \eqref{eq_identity_u1}}, we can write the local state cost functions \eqref{eq_state_cost_primal_dual} \lhe{for each $i \in \calv$ and $j \in \cale$} in the form
\begin{align*}
    \begin{split}
        q_i(\xi_i,\tau_i,R_i) & = \left(\theta_i - \theta_i^{\ast} \right) \left(\nabla F_i(\theta_i) - \nabla F_i(\theta_i^{\ast}) + \eta_i - \eta_i^{\ast} \right)  \\
            & \qquad + \textstyle \sum_{k = 2}^{\rho_i} \frac{a_{ik}}{2b_{ik}} \xi_{ik}^2 \\
        & = \nabla^T V_i(\xi_i) b_i \mathbb{1}_{\rho_i} \left[ \nabla F_i(\theta_i) - \nabla F_i(\theta_i^{\ast}) \right. \\
            & \qquad \left. + \eta_i - \eta_i^{\ast} \right]
            +\lVert \underline{u}_i \rVert^2_{R_i},
    \end{split} \\
    \nonumber \\
    \begin{split}
        \hat{q}_i(\tau_i,\xi_i,\hat{R}_i)  &= -(\tau_{i1} - \lambda_i^{\ast})[G_i(\theta_i)]_{\tau_{i1}}^+ \\
            & \qquad - \textstyle \sum_{k = 2}^{\hat{\rho}_i} \tau_{ik}G_i(\theta_i) + \textstyle \sum_{k = 2}^{\hat{\rho}_i} \frac{\hat{a}_{ik}}{2\hat{b}_{ik}} \tau_{ik}^2 \\
        & = -(\nabla \hat{V}_i(\tau_i))_1 [\hat{b}_{i1} G_i(\theta_i)]_{\tau_{i1}}^+ \\
            & \qquad - \textstyle  \sum_{k = 2}^{\hat{\rho}_i} (\nabla \hat{V}_i(\tau_i))_k \hat{b}_{ik} G_i(\theta_i)
        + \lVert \underline{\hat{u}}_i \rVert_{\hat{R}_i}^2 ,
    \end{split} \\
    \nonumber \\
    \begin{split}
        \check{q}_j(\zeta_j,\check{R}_j)  &= \textstyle \sum_{k =2}^{\check{\rho}_j} \frac{\check{a}_{jk}}{2\check{b}_{jk}} \zeta_{jk}^2 =
        \lVert \underline{\check{u}}_j \rVert_{\check{R}_j}^2 .
    \end{split}
\end{align*}
Therefore, \lhe{using the expressions above in \eqref{eq_full_state_cost} gives},
\begin{equation*}
    \begin{split}
        q(\xi,\tau,\zeta, & R)
        = \textstyle \sum_{i \in \calv} \left[  -\nabla^T V_i(\xi_i) b_i \mathbb{1}_{\rho_i} \left(-\nabla F_i(\theta_i) \right. \right. \\
            & \qquad \qquad  \left. \left. + \nabla F_i(\theta_i^\ast) - \eta_i + \eta_i^\ast
            \right) + \lVert \underline{u}_i \rVert^2_{R_i} \right]
        \\
        & \qquad + \textstyle \sum_{i \in \calv} \left[ -(\nabla \hat{V}_i(\tau_i))_1 [\hat{b}_{i1} G_i(\theta_i)]_{\tau_{i1}}^+ \right. \\
            & \qquad \qquad   \left.- \textstyle  \sum_{k = 2}^{\hat{\rho}_i} (\nabla \hat{V}_i(\tau_i))_k \hat{b}_{ik} G_i(\theta_i)
            +\lVert \underline{\hat{u}}_i \rVert_{\hat{R}_i}^2  \right]
        \\ & \qquad+ \textstyle \sum_{j \in \cale}
        \lVert \underline{\check{u}}_j \rVert^2_{\check{R}_j} .
    \end{split}
\end{equation*}
\lhe{Next, by adding \eqref{eq_proof_relation} we get
\begin{equation*}
    \begin{split}
        q(\xi,\tau, &\zeta,  R)
        = \textstyle \sum_{i \in \calv} \left[ -\nabla^T V_i(\xi_i) b_i \mathbb{1}_{\rho_i} \left(\phi_i + \nabla F_i(\theta_i^\ast) \right. \right.  \\
            & \qquad \qquad \left. \left. + \eta_i^\ast  + \psi^\ast
            \right) + \lVert \underline{u}_i \rVert^2_{R_i} \right]
        \\
        & \qquad + \textstyle \sum_{i \in \calv} \left[ -(\nabla \hat{V}_i(\tau_i))_1 [\hat{b}_{i1} G_i(\theta_i)]_{\tau_{i1}}^+ \right. \\
            &  \qquad \qquad \left. - \textstyle  \sum_{k = 2}^{\hat{\rho}_i} (\nabla \hat{V}_i(\tau_i))_k \hat{b}_{ik} G_i(\theta_i)
            +\lVert \underline{\hat{u}}_i \rVert_{\hat{R}_i}^2  \right]
        \\
        & \qquad + \textstyle \sum_{j \in \cale} \left[-\nabla^T \check{V}_j(\zeta_j) \check{b}_j \mathbb{1}_{\check{\rho}_j} \left(\omega_j - \omega_j^\ast \right) \right. \\
        & \qquad \qquad \left. +\lVert \underline{\check{u}}_j \rVert^2_{\check{R}_j} \right].
    \end{split}
\end{equation*}
\lhk{Finally, applying the KKT conditions \eqref{eq_distributed_kkt} gives}
\begin{equation*} 
    \begin{split}
        q(\xi,\tau,\zeta, & R)
        = \textstyle \sum_{i \in \calv} \left[ -\nabla^T V_i(\xi_i)  b_i \mathbb{1}_{\rho_i} \phi_i
        +\lVert \underline{u}_i \rVert^2_{R_i} \right]
        \\
        & \qquad + \textstyle \sum_{i \in \calv} \left[ -(\nabla \hat{V}_i(\tau_i))_1 [\hat{b}_{i1} G_i(\theta_i)]_{\tau_{i1}}^+ \right. \\
        & \qquad \left. - \textstyle  \sum_{k = 2}^{\hat{\rho}_i} (\nabla \hat{V}_i(\tau_i))_k \hat{b}_{ik} G_i(\theta_i)
        + \lVert \underline{\hat{u}}_i \rVert_{\hat{R}_i}^2  \right]
        \\
        & \qquad + \textstyle \sum_{j \in \cale} \left[-\nabla^T \check{V}_j(\zeta_j)  \check{b}_j \mathbb{1}_{\check{\rho}_j} \omega_j
        +\lVert \underline{\check{u}}_j \rVert^2_{\check{R}_j} \right],
    \end{split}
\end{equation*}
where $\phi_i$ is defined in \eqref{eq_def_phi}.}
\\ \hbox{} \\
\textbf{Step 4:}
\lhe{Next, we evaluate the term $\lVert U \rVert_R^2$ in \eqref{eq_J} with $U$ in \eqref{eq_define_v} and $R$ as in \eqref{eq_control_cost_ineq}:
\begin{equation*} 
    \begin{split}
        \lVert U \rVert_R^2
        & = \lVert \underline{U} \rVert_R^2 + 2 (\underline{U})^T R v + \lVert v \rVert_R^2 \\
        &= 2(\underline{U})^T R U - \lVert \underline{U} \rVert_R^2  + \lVert v \rVert_R^2 \\
        & = \textstyle \sum_{i \in \calv} \left( 2 \underline{u}_i^T R_i u_i - \lVert \underline{u}_i \rVert_{R_i}^2 \right) \\
            & \qquad + \textstyle \sum_{i \in \calv} \left( 2 \underline{\hat{u}}_i^T \hat{R}_i \hat{u}_i - \lVert \underline{\hat{u}}_i \rVert_{\hat{R}_i}^2 \right) \\
            & \qquad
             + \textstyle \sum_{j \in \cale} \left( 2 \underline{\check{u}}_j^T \check{R}_j \check{u}_j - \lVert \underline{\check{u}}_j \rVert_{\check{R}_j}^2 \right)
            + \lVert v \rVert_R^2 \\
        & = \textstyle \sum_{i \in \calv} \left( -\nabla^T V_i(\xi_i) b_{ui} u_i - \lVert \underline{u}_i \rVert_{R_i}^2 \right) \\
            & \qquad + \textstyle \sum_{i \in \calv} \left( -\nabla^T \hat{V}_i(\tau_i) \hat{b}_{ui} \hat{u}_i - \lVert \underline{\hat{u}}_i \rVert_{\hat{R}_i}^2 \right)
            \\ & \qquad
            + \textstyle \sum_{j \in \cale} \left(-\nabla^T \check{V}_j(\zeta_j) \check{b}_{uj} \check{u}_j - \lVert \underline{\check{u}}_j \rVert_{\check{R}_j}^2 \right)
           \! + \! \lVert v \rVert_R^2,
    \end{split}
\end{equation*}
where the identities \eqref{eq_identity_u0} were used in the last equality.}
\\ \hbox{} \\
\textbf{Step 5:} \lhe{We are now ready to evaluate $J(U)$ in \eqref{eq_J} using $U$ in \eqref{eq_define_v}.}
\lhe{Combining the expressions for $q(\xi,\tau,\zeta,R)$ and $\lVert U \rVert_R^2$ gives:}
\begin{align*}
        J(U) & = \int_0^\infty \lVert U \rVert_R^2 + q(\xi,\tau,\zeta,R) \ dt \\
        & = - \sum_{i \in \calv} \int_0^\infty  \nabla^T V_i(\xi_i) (b_i \mathbb{1}_{\rho_i} \phi_i + b_{ui} u_i) \ dt
        \\ & \qquad
        -\sum_{i \in \calv} \int_0^\infty \left[ (\nabla \hat{V}_i(\tau_i))_1[\hat{b}_{i1}G_i(\theta_i)]_{\tau_{i1}} \right. \\
        & \qquad \qquad \left. + \textstyle \sum_{k = 2}^{\hat{\rho}_i} (\nabla \hat{V}_i(\tau_i))_k (\hat{b}_{ik}G_i(\theta_i)+ \hat{u}_{ik}) \right] \ dt \\
        & \qquad - \sum_{j \in \cale} \int_0^\infty    \nabla^T \check{V}_j(\zeta_j) (\check{b}_j \mathbb{1}_{\check{\rho}_j} \omega_j+ \check{b}_{uj} \check{u}_j) \ dt
        \\& \qquad
        + \int_0^\infty \lVert v \rVert_R^2  \ dt.
\end{align*}
\lhh{Now,} for each $i \in \calv$ and for $k \in \{2, \ldots, \hat{\rho}_i\}$
\begin{equation} \label{eq_proj_equality}
    \begin{split}
        (\nabla \hat{V}_i (\tau_i) )_k (& \hat{b}_{ik} G_i(\theta_i) +\hat{u}_{ik}) \\ & = ( \nabla \hat{V}_i(\tau_i) )_k \left[\hat{b}_{ik}G_i(\theta_i)+\hat{u}_{ik} \right]_{\tau_i}^+ ,
    \end{split}
\end{equation}
where $( \nabla \hat{V}_i(\tau_i) )_k$ is the $k^{\mathrm{th}}$ component of $\nabla \hat{V}_i(\tau_i)$.
To see this, firstly note that $( \nabla \hat{V}_i(\tau_i) )_k = \tfrac{\tau_{ik}}{\hat{b}_{ik}}$ for $k \in \{2, \ldots, \hat{\rho}_i\}$. Using \eqref{eq_positive_projection}, equality holds for $\tau_{ik} > 0$. For $\tau_{ik} = 0$, we have $( \nabla \hat{V}_i(\tau_i) )_k = 0$, meaning equality also holds in all other cases.
\lhh{Therefore}, using \eqref{eq_proj_equality}, we can write the terms in the second integrand in 
\lhk{$J(U)$} for each $i \in \calv$ as
\begin{align*}
    & (\nabla \hat{V}_i(\tau_i))_1[\hat{b}_{i1}G_i(\theta_i)]_{\tau_{i1}}^+ \\
        & \qquad + \textstyle \sum_{k = 2}^{\hat{\rho}_i} (\nabla \hat{V}_i(\tau_i))_k (\hat{b}_{ik}G_i(\theta_i)+ \hat{u}_{ik}) \\
    &  = (\nabla \hat{V}_i(\tau_i))_1[\hat{b}_{i1}G_i(\theta_i)]_{\tau_{i1}}^+ \\
    & \qquad + \textstyle \sum_{k = 2}^{\hat{\rho}_i} (\nabla \hat{V}_i(\tau_i))_k [\hat{b}_{ik}G_i(\theta_i )+ \hat{u}_{ik}]_{\tau_{ik}}^+ \\
    & = \nabla^T \hat{V}_i(\tau_i) [\hat{b}_{i}\mathbb{1}_{\hat{\rho}_i}G_i(\theta_i)+ \hat{b}_{ui}\hat{u}_{i}]_{\tau_{i}}^+ \\
    &= \nabla^T \hat{V}_i(\tau_i) \dot{\tau}_i,
\end{align*}
\lhe{where here we substituted the dynamics \licc{in} \eqref{eq_primal_dual_control_lambda}. Similarly, substituting the dynamics \licc{in} \eqref{eq_primal_dual_control_theta} and \eqref{eq_primal_dual_control_mu} into $J(U)$ gives}
\begin{align*}
        J(U)
        & = \int_0^\infty \left[ - \sum_{i \in \calv} \nabla^T V_i(\xi_i) \dot{\xi}_i
        - \sum_{i \in \calv}  \nabla^T \hat{V}_i(\tau) \dot{\tau}_i \right. \\
            & \qquad \qquad \left.- \sum_{j \in \cale}  \nabla \check{V}_j(\zeta_j) \dot{\zeta}_j \right] \ dt
            + \int_0^\infty v^T R v \ dt \\
        & = \int_0^\infty \left[- \sum_{i \in \calv}    \dfrac{dV_i(\xi_i)}{dt}
        - \sum_{i \in \calv}  \frac{d\hat{V}_i(\tau_i)}{dt} \right. \\
            & \qquad \qquad \left. - \sum_{j \in \cale}   \frac{d\check{V}_j(\zeta_j)}{dt} \right] \ dt
            + \int_0^\infty v^T R v \ dt \\
        & = V(\xi_0,\tau_0,\zeta_0) + \int_0^\infty v^T R v \ dt,
\end{align*}
where here we evaluated $ \lim_{t \rightarrow \infty} V(\xi,\tau,\zeta) = 0$ as we assume that the control action $U$ is \lhj{stabilizing}. Therefore, we have that \lhe{$J(U)$ with \eqref{eq_define_v} evaluates to a constant term (\lhj{dependent} on the initial conditions $(\xi_0,\tau_0,\zeta_0)$) plus an integral that depends on the arbitrary vector $v$. We see that this term} is \lhj{minimized} for $v = 0$ \lhe{(as $R$ is positive definite), meaning that $J(U)$ is \lhj{minimized} by}
$U = \underline{U}$ \lha{given by \eqref{eq_augmented_controllers}}.

\section{Proof of Lemmas \ref{lemma_state_cost_bound} and \ref{lemma_non_neg}} \label{sec_proof_lemmas}
\subsection{Proof of Lemma \ref{lemma_state_cost_bound}}
    \lhb{We see \eqref{eq_term_interp1} follows directly from \eqref{eq_bregman_def}. To show \eqref{eq_term_interp2}, we first note}
    \begin{equation*}
        (\tau_{i1} - \lambda_i^{\ast}) [G_i(\theta_i)]_{\tau_{i1}}^+ \leq (\tau_{i1} - \lambda_i^{\ast}) G_i(\theta_i),
    \end{equation*}
    which follows directly from the definition of the operator $[\cdot]_{\ast}^+$ in \eqref{eq_positive_projection} (equality holds except in the case $\tau_{i1}=0$ and $G_i(\theta_i) <0$, which means the left-hand side evaluates to zero, while the right-hand side evaluates to a non-negative number). Therefore,
    \begin{equation*}
        \begin{split}
            \sigma_{i4} & \geq - (\tau_{i1} - \lambda_i^{\ast}) G_i(\theta_i) - \textstyle \sum_{k = 2}^{\hat{\rho}_i} \tau_{ik}G_i(\theta_i) \\
            & = - (\lambda_i - \lambda_i^{\ast})G_i(\theta_i).
        \end{split}
    \end{equation*}
    Now, \lhb{
    \begin{equation*}
        \begin{split}
            \sigma_{i2} + \sigma_{i4} &\geq (\theta_i - \theta_i^{\ast})(\eta_i - \eta_i^{\ast})  - (\lambda_i - \lambda_i^{\ast}) G_i(\theta_i) \\
            & = \left[(\theta_i - \theta_i^{\ast}) \nabla G_i(\theta_i) \right] \lambda_i
                \\ &\qquad -\left[(\theta_i \!-\! \theta_i^{\ast})\nabla G_i(\theta_i^{\ast})\right]\lambda_i^{\ast}
             - (\lambda_i \!-\! \lambda_i^{\ast})G_i(\theta_i) \\
            & = \left[ D_{G_i}(\theta_i^\ast,\theta_i) -G_i(\theta_i^\ast)+G_i(\theta_i)\right] \lambda_i
                \\ &\qquad + \left[ D_{G_i}(\theta_i,\theta_i^\ast) - G_i(\theta_i) + G_i(\theta_i^\ast)\right] \lambda_i^\ast\\
                & \qquad - (\lambda_i - \lambda_i^{\ast})G_i(\theta_i) \\
            & = D_{G_i} (\theta_i^\ast,\theta_i)\lambda_i + D_{G_i} (\theta_i,\theta_i^{\ast}) \lambda_i^\ast  - G_i(\theta_i^\ast) \lambda_i \\
        \end{split}
    \end{equation*}}%
    where we have used
    $G_i(\theta_i^{\ast})\lambda_i^{\ast} = 0$ \lhb{from \eqref{eq_distributed_kkt_2} and
    the definition of the Bregman divergence \eqref{eq_bregman_def}, which gives us \eqref{eq_term_interp2}.}

\subsection{Proof of Lemma \ref{lemma_non_neg}}
    \lhb{Firstly, we note that the term $\lVert U \rVert_R^2 \lic{\geq} 0$ 
    as $R$ is positive definite. Next, we see that $\sigma_{i1} \geq 0$ via \eqref{eq_term_interp1}, as the Bregman divergence is non-negative by the convexity of $F_i{(\cdot)}$. Similarly, $\sigma_{i2}+\sigma_{i4}$ in \eqref{eq_term_interp2} is non-negative, as the Bregman divergence terms are non-negative,}
    $\lambda_i \geq 0$ and $G_i(\theta_i^{\ast}) \leq 0$.
    \lhb{Finally, the terms $\sigma_{i3}$, $\sigma_{i5}$ and $\sigma_{j1}$ in \eqref{eq_state_cost_primal_dual} are all non-negative. Therefore, the cost in \eqref{eq_primal_cost} is non-negative.}

    \lhb{We note that \eqref{eq_primal_dual_control} with \eqref{eq_augmented_controllers} recovers the algorithm \eqref{eq_comp_dynamics}. Equilibrium for \eqref{eq_comp_dynamics} is reached when, for each $i \in \calv$ and $j \in \cale$, $\xi_{ik}=0$, $\tau_{ik}=0$ and $\zeta_{jk}=0, \forall k \geq 2$, and $\xi_{i1} = \theta_i^\ast$, $\tau_{i1}=\lambda_i^\ast$ and $\zeta_{j1} = \mu_j^\ast$. Therefore, at equilibrium, $U^\ast = 0$, and each term in \eqref{eq_state_cost_primal_dual} evaluates to zero, meaning the cost function in \eqref{eq_primal_cost} is zero.}

\section{Proof of Proposition \ref{prop_feed_forward}} \label{sec_proof_ff}
The proof follows that of Proposition \ref{prop_main_result}.
\lhi{As discussed in Section \ref{sec_feedforward}, for simplicity \lic{in} the presentation, we consider the unconstrained case of \eqref{opt_distributed_transformed} as in Corollary \ref{cor_unconstrained_no_ff}, i.e., with no local inequality constraints, so $\Sigma_\lambda^+=0$ \lic{in Figure \ref{fig_dist_opt}}. With this simplification, the dual dynamics associated with the local inequality constraints will be omitted.}
We note that the problem including \lhi{local} inequality constraints can be considered through an analogous optimal control problem, which contains state cost functions similar to \eqref{eq_state_cost_primal_dual}.

\subsection{Preliminary Definitions}
We first introduce some notation that will be useful within the proof. The control cost matrix is given by
\begin{align}
    R  = \mathrm{diag}\begin{bmatrix}
        R_\calv, &  R_\cale
    \end{bmatrix}, \label{eq_R_ff}
\end{align}
where $R_\calv$ and $R_\cale$ are given in \eqref{eq_opt_ext_R}. Next, for each $i \in \calv$, let $\vec{\theta}_i^\ast = [\theta_i^\ast, \mathbb{0}_{\rho_i-1}^T]^T$, and for each $j \in \cale$, let $\vec{\mu}_j^\ast = [\mu_j^\ast,\mathbb{0}_{\check{\rho}_j-1}^T]^T$. Then let $\vec{\theta}^\ast = [\vec{\theta}_i^\ast]_{i\in\calv}$ and $\vec{\mu}^\ast = [\vec{\mu}_j^\ast]_{j\in\cale}$. With this notation, we take the aggregate storage functions associated with the node and edge subsystems respectively as
\begin{subequations}
    \begin{align}
        V_\calv(\xi) & = \sum_{i \in \calv} V_i(\xi_i) = \frac{1}{2} (\xi-\vec{\theta}^\ast)^T b^{-1} (\xi-\vec{\theta}^\ast) \\
        V_\cale(\zeta) & = \sum_{j \in \cale} V_j(\zeta_j) = \frac{1}{2} (\zeta - \vec{\mu}^\ast)^T \check{b}^{-1}(\zeta - \vec{\mu}^\ast),
    \end{align}
\end{subequations}
where $V_i(\xi_i)$ and $V_j(\zeta_j)$ are given in \eqref{eq_storage_local_1} and \eqref{eq_storage_local_3}.
Finally, the cost functional \licc{in} \eqref{eq_opt_ext_cost} to be \lhj{minimized} is defined as
\begin{equation} \label{eq_J_ff}
    \begin{split}
         J(U): = \int_0^{\infty} \lVert U \rVert_R^2 + q(\xi,\zeta,R) \ dt,
    \end{split}
\end{equation}
where $U = \begin{bmatrix} u_\calv^T & u_\cale^T \end{bmatrix}^T$, with $u_\calv$ and $u_\cale$ the control inputs for dynamics \eqref{eq_primal_dual_control_uncon_ff}, $R$ is defined in \eqref{eq_R_ff}, and
\begin{equation} \label{eq_full_state_cost_ff}
    q(\xi,\zeta,R) = q_\calv(\xi,R_\calv) + q_\cale(\zeta,R_\cale),
\end{equation}
with $q_\calv(\xi,R_\calv)$ and $q_\cale(\zeta,R_\cale)$ defined in \eqref{eq_state_cost_ff}.

\subsection{Proof}
As in the proof of Proposition \ref{prop_main_result}, we directly evaluate the cost functional $J(U)$ in \eqref{eq_J_ff} under the \lhj{stabilizing} control action given by
\begin{equation} \label{eq_def_v_ff}
    U = \underline{U} + v,
\end{equation}
where $\underline{U} = \begin{bmatrix} \underline{u}_\calv^T & \underline{u}_\cale^T \end{bmatrix}^T$ is our proposed optimal controller with $\underline{u}_\calv$ and $\underline{u}_\cale$ given in \eqref{eq_opt_ext_controller} and $v$ is an \ils{arbitrary function of time} of appropriate dimension.

We now replicate the steps in Section \ref{sec_proof_p1_steps}, updating the notation and cost functional to align with those of Proposition \ref{prop_feed_forward}.
\\ \hbox{} \\
\textbf{Step 1:} We repeat the analysis of Step 1 in Section \ref{sec_proof_p1_steps}. Using the notation for the aggregate node and edge subsystems (defined above and in Section \ref{sec_feedforward}), we can write
    \begin{align}
        \nabla V_\calv(\xi) &= b^{-1} (\xi-\vec{\theta}^\ast), &
        \nabla V_\cale(\zeta) & = \check{b}^{-1}(\zeta - \vec{\mu}^\ast),
    \end{align}
which gives
    \begin{align} \label{eq_storage_gradient_ff}
        \nabla^T V_\calv(\xi) b \mathbb{1}_\calv & = (\theta-\theta^{\ast})^T,  &
        \nabla^T V_\cale(\zeta) \check{b} \mathbb{1}_\cale & = (\tilde{\mu}-\mu^{\ast})^T.
    \end{align}
Then, similarly to \eqref{eq_proof_relation}, we get
\begin{equation} \label{eq_proof_relation_ff}
     \nabla^T V_\calv(\xi) b \mathbb{1}_\calv (\tilde{\psi} - \psi^{\ast}) \\
     - \nabla^T V_\cale(\zeta) \check{b} \mathbb{1}_\cale (\omega - \omega^{\ast}) = 0.
\end{equation}
\textbf{Step 2:} Now, we obtain identities relating the proposed optimal control to our storage functions, using notation from above and Section \ref{sec_feedforward}. Firstly, for the node subsystem, we have:
\begin{equation} \label{eq_ff_u_identity_1}
    \begin{split}
            -\frac{1}{2} &  R_\calv^{-1} b_{u \calv}^T \nabla V_\calv(\xi) \\
            & = - \frac{1}{2} \begin{bmatrix} 2a  I_{o\calv}^T b^{-1} I_{o \calv} & 0 \\ 0 & 2 \check{d} \end{bmatrix} \begin{bmatrix} I_{o\calv}^T \\ \cala^T \mathbb{1}_\calv \end{bmatrix} b b^{-1} \left(\xi - \vec{\theta}^{\ast} \right) \\
            & = - \begin{bmatrix} a I_{o \calv}^T b^{-1} \\ \check{d} \cala^T \mathbb{1}_\calv^T \end{bmatrix} \xi  = \underline{u}_\calv,
    \end{split}
\end{equation}
where here we used identities
\begin{align}
    & I_{o\calv}^T \vec{\theta}^\ast = \left( \oplus_{i\in\calv} I_{o\rho_i} \right) [\vec{\theta}_i^\ast]_{i\in\calv} = [\mathbb{0}_{\rho_i}]_{i\in\calv}, \label{eq_ff_u_identity_2} \\
    & \begin{aligned}
    \cala^T \mathbb{1}_\calv^T \vec{\theta}^\ast
    & = \cala^T \left( \oplus_{i\in\calv} \mathbb{1}_{\rho_i}^T \right) [\vec{\theta}_i^\ast]_{i\in\calv} \\
    & = \cala^T [\theta_i^\ast]_{i\in\calv} = \cala^T \theta^\ast = 0, \label{eq_ff_u_identity_3}
    \end{aligned}
\end{align}
where the last statement follows via the KKT condition \eqref{eq_distributed_kkt_3}. Similarly, for the edge subsystem, we have
\begin{equation}  \label{eq_ff_u_identity_4}
    \begin{split}
        -\frac{1}{2} & R_\cale^{-1} b_{u \cale}^T \nabla^T V_\cale(\zeta) \\
             & = - \frac{1}{2} \left( 2 \check{a}  I_{o\cale}^T \check{b}^{-1} I_{o \cale} \right) \left( I_{o\cale}^T \check{b} \right) \check{b}^{-1} \left(\zeta - \vec{\mu}^{\ast} \right)
             \\ & = - \check{a} I_{o \cale}^T \check{b}^{-1} \zeta  = \underline{u}_\cale,
    \end{split}
\end{equation}
where we used identities
\begin{align}
    &I_{o\cale}^T \vec{\mu}^\ast = \left( \oplus_{j\in\cale} I_{o\rho_j} \right) [\vec{\mu}_j^\ast]_{j\in\cale} = [\mathbb{0}_{\rho_j}]_{j\in\cale}, \label{eq_ff_u_identity_5} \\
    & I_{o \cale}^T \check{b}^{-1} I_{o \cale} I_{o \cale}^T = I_{o \cale}^T \check{b}^{-1}. \label{eq_ff_u_identity_6}
\end{align}
Therefore, using \eqref{eq_ff_u_identity_1}, we have
 \begin{align*}
    \label{eq_ff_u_identity_7}
        \lVert \underline{u}_\calv \rVert_{R_\calv}^2
        & = \frac{1}{4} \nabla^T V_\calv(\xi) b_{u\calv} R_\calv^{-1} R_\calv  R_\calv^{-1} b_{u \calv}^T \nabla V_\calv(\xi) \\
        & = - \frac{1}{2} \nabla^T V_\calv(\xi) b_{u\calv} \underline{u}_\calv \\
        & = \frac{1}{2} (\xi-\vec{\theta}^\ast)^T b^{-1} b \begin{bmatrix} I_{o\calv} & \mathbb{1}_\calv \cala \end{bmatrix} \begin{bmatrix}  a I_{o \calv}^T b^{-1} \\ \check{d} \cala^T \mathbb{1}_\calv^T  \end{bmatrix} \xi \\
        &=  \frac{1}{2} \xi^T \mathbb{1}_\calv \cala \check{d}\cala^T \mathbb{1}_\calv^T \xi + \frac{1}{2} \xi^T I_{o\calv} a I_{o\calv}^T b^{-1} \xi \\
        & = \frac{1}{2} \theta^T \mathcal{L}_{\check{d}} \theta + \sum_{i\in\calv} \sum_{k=2}^{\rho_i} \frac{a_{ik}}{2b_{ik}} \xi_{ik}^2,
 \end{align*}
where here we used \eqref{eq_ff_u_identity_2}, \eqref{eq_ff_u_identity_3}, and $\mathbb{1}_\calv^T \xi = \left( \oplus_{i \in \calv} \mathbb{1}_{\rho_i}^T \right) [\xi_i]_{i\in\calv} = [\mathbb{1}_{\rho_i}^T \xi_i]_{i\in\calv} = [\theta_i]_{i\in\calv} = \theta $.
\lhk{Following a similar approach using \eqref{eq_ff_u_identity_4} gives
\begin{equation} \label{eq_ff_u_identity_8}
    \lVert \underline{u}_\cale \rVert_{R_\cale}^2 = \sum_{j\in\cale} \sum_{k=2}^{\check{\rho}_j} \frac{\check{a}_{jk}}{2\check{b}_{jk}} \zeta_{jk}^2.
\end{equation}
}%
\textbf{Step 3:} We now write the state cost functions \eqref{eq_state_cost_ff} in terms of our storage functions and proposed optimal controllers:
    \begin{align*}
       \begin{split}
            q_\calv(\xi,R_\calv)
            & = \textstyle \frac{1}{2}\theta^T \mathcal{L}_{\check{d}} \theta +   \sum_{i \in \calv} \left[ \left(\theta_i - \theta_i^{\ast} \right) \left(\nabla F_i(\theta_i) \right. \right. \\
                & \left. \left. \qquad \textstyle - \nabla F_i(\theta_i^{\ast}) \right)  + \sum_{k = 2}^{\rho_i} \frac{a_{ik}}{2b_{ik}} \xi_{ik}^2 \right] \\
            & = (\theta - \theta^\ast)^T(\nabla F(\theta) - \nabla F(\theta^\ast)) \\
                & \qquad \textstyle + \frac{1}{2}\theta^T \mathcal{L}_{\check{d}} \theta + \sum_{i\in\calv} \sum_{k = 2}^{\rho_i} \frac{a_{ik}}{2b_{ik}} \xi_{ik}^2 \\
            & = \nabla^T V_\calv(\xi) b \mathbb{1}_\calv (\nabla F(\theta) - \nabla F(\theta^\ast)) + \lVert \underline{u}_\calv \rVert_{R_\cale}^2,
        \end{split}   \\
        \nonumber \\
        \begin{split}
             q_\cale(\zeta,R_\cale)
            &= \sum_{j \in \cale} \check{q}_j(\zeta_j,\check{R}_j) \\
            & = \sum_{j\in\cale} \sum_{k =2}^{\check{\rho}_j} \frac{\check{a}_{jk}}{2\check{b}_{jk}} \zeta_{jk}^2 \\
            & = \lVert \underline{u}_\cale \rVert_{R_\cale}^2.
        \end{split}
    \end{align*}
Therefore, by adding \eqref{eq_proof_relation_ff} to \eqref{eq_full_state_cost_ff}, we get
\begin{equation*} 
    \begin{split}
        q(\xi,\zeta,R)
        & = \nabla^T V_\calv(\xi) b \mathbb{1}_\calv (\nabla F(\theta) - \nabla F(\theta^\ast) + \tilde{\psi}-\psi^\ast)  \\
            & \quad  + \lVert \underline{u}_\calv \rVert_{R_\cale}^2 - \nabla^T V_\cale(\zeta) \check{b} \mathbb{1}_\cale (\omega - \omega^{\ast}) + \lVert \underline{u}_\cale \rVert_{R_\cale}^2 \\
        & = - \nabla^T V_\calv(\xi) b \mathbb{1}_\calv \tilde{\phi}  + \lVert \underline{u}_\calv \rVert_{R_\cale}^2 \\
            & \quad  - \nabla^T V_\cale(\zeta) \check{b} \mathbb{1}_\cale \omega + \lVert \underline{u}_\cale \rVert_{R_\cale}^2,
    \end{split}
\end{equation*}
where $\tilde{\phi}$ is given in \eqref{eq_def_psi_phi_tilde}, and we used the KKT conditions
        $ \nabla F(\theta^\ast) + \psi^\ast = 0$ and
        $ \omega^\ast = 0$.
\\ 
\textbf{Step 4:} We can evaluate the term $\lVert U \rVert_R^2$ in \eqref{eq_J_ff} as
\begin{equation*} 
    \begin{split}
        \lVert U \rVert_R^2
        & = 2(\underline{U})^T R U - \lVert \underline{U} \rVert_R^2  + \lVert v \rVert_R^2 \\
        & = 2 \underline{u}_\calv^T R_\calv u_\calv  + 2 \underline{u}_\cale^T R_\cale u_\cale \\
            & \qquad\qquad - \lVert \underline{u}_\calv \rVert_{R_\calv} - \lVert \underline{u}_\cale \rVert_{R_\cale} + \lVert v \rVert_R^2 \\
        & = -\nabla^T V_\calv(\xi) b_{u\calv} u_\calv - \nabla^T V_\cale(\zeta) b_{u\cale} u_\cale \\
            & \qquad\qquad- \lVert \underline{u}_\calv \rVert_{R_\calv} - \lVert \underline{u}_\cale \rVert_{R_\cale} + \lVert v \rVert_R^2
    \end{split}
\end{equation*}
\\ \hbox{} \\
\textbf{Step 5:} Therefore, by combining the expressions for $q(\xi,\zeta,R)$ and $\lVert U \rVert_R^2$, we can calculate cost functional \eqref{eq_J_ff} with controller \eqref{eq_def_v_ff} as
\begin{align*}
        J(U)
        & = \int_0^\infty \left[-\nabla^T V_\calv(\xi) \left( b \mathbb{1}_\calv \tilde{\phi} + b_{u\calv} u_\calv \right) \right. \\
            & \qquad \qquad \left. - \nabla^T V_\cale(\zeta) \left( \check{b} \mathbb{1}_\cale \omega + b_{u\cale} u_\cale \right) \right] \ dt \\
             & \qquad \qquad + \int_0^\infty \lVert v \rVert_R^2 \ dt \\
        & = - \int_0^\infty \frac{dV_\calv(\xi)}{dt} \ dt -  \int_0^\infty \frac{dV_\cale(\zeta)}{dt} \ dt \\
            & \qquad \qquad + \int_0^\infty \lVert v \rVert_R^2 \ dt \\
        & = V_\calv(\xi_0) + V_\cale(\zeta_0) + \int_0^\infty \lVert v \rVert_R^2 \ dt,
\end{align*}
where here we evaluated $\lim_{t \rightarrow \infty} V_\calv(\xi) = 0$ and $\lim_{t \rightarrow \infty} V_\cale(\zeta) = 0$ as we assume that the control action $U$ is \lhj{stabilizing}. Therefore, we have that $J(U)$ with \eqref{eq_def_v_ff} evaluates to a constant term (\lhj{dependent} on the initial conditions $(\xi_0,\zeta_0)$) plus an integral that depends on the arbitrary vector $v$. We see that this term is \lhj{minimized} for $v = 0$ (as $R$ is positive definite), meaning that $J(U)$ is \lhj{minimized} by
$U = \underline{U}$ \lha{given by \eqref{eq_opt_ext_controller}}.

\begin{ack}                               
This work was funded by UKRI grant EP/T517847/1. For the purpose of open access, the authors have applied a Creative Commons Attribution (CC BY) licence to any Author Accepted Manuscript version arising.
\end{ack}

\bibliographystyle{plain}        
\bibliography{Files/bib_auto_upload}           

\end{document}

%% file: Files/preamble.tex
\newcommand{\calv}{\mathcal{V}}
\newcommand{\cale}{\mathcal{E}}

\newcommand{\cala}{\mathcal{A}}

\newcommand{\il}[1]{\textcolor{black}{#1}} 
\newcommand{\ic}[1]{\textcolor{black}{#1}} 
\newcommand{\lhc}[1]{\textcolor{black}{#1}} 
\newcommand{\lhm}[1]{\textcolor{black}{#1}} 
\newcommand{\icl}[1]{\textcolor{black}{#1}} 
\newcommand{\lhg}[1]{\textcolor{black}{#1}} 
\newcommand{\icc}[1]{\textcolor{black}{#1}} 
\newcommand{\lhf}[1]{\textcolor{black}{#1}} 
\newcommand{\iccc}[1]{\textcolor{black}{#1}} 
\newcommand{\lhab}[1]{\textcolor{black}{#1}} 
\newcommand{\lhad}[1]{\textcolor{black}{#1}} 

\newcommand{\lha}[1]{\textcolor{black}{#1}} 
\newcommand{\lhb}[1]{\textcolor{black}{#1}} 
\newcommand{\lhe}[1]{\textcolor{black}{#1}} 
\newcommand{\lhh}[1]{\textcolor{black}{#1}} 
\newcommand{\lhi}[1]{\textcolor{black}{#1}} 
\newcommand{\lhj}[1]{\textcolor{black}{#1}} 

\newcommand{\li}[1]{\textcolor{black}{#1}} 
\newcommand{\lic}[1]{\textcolor{black}{#1}} 
\newcommand{\licc}[1]{\textcolor{black}{#1}} 

\usepackage[prependcaption,colorinlistoftodos]{todonotes}

\usepackage{amsfonts}
\usepackage{amssymb}
\usepackage{graphicx}
\usepackage{subcaption}
\usepackage{epstopdf}
\usepackage{algorithmic}
\usepackage{bbold}
\usepackage[short]{optidef}
\usepackage{tikz}

\usepackage[sort]{cite}

\newtheorem{proposition}[thm]{Proposition}
\newtheorem{corollary}[thm]{Corollary}
\newtheorem{lemma}[thm]{Lemma}

\newtheorem{remark}{Remark}

%% file: Images/tikz_block_diagram.tex
\begin{tikzpicture}[transform shape,scale=1]


\node[matrix,column sep=0.5cm, row sep =0.25cm] at (3.75,0)
{
    \node[draw, circle, minimum size=0.05cm](sum_1){};		
			\node at (sum_1.north west)[left=-2pt,yshift=+3pt]{ $-$};
			\node at (sum_1.south east)[right=-2pt,yshift=-3pt]{$-$}; &
    \node[draw, rectangle,
			minimum width=0.75cm, 	minimum height=0.75cm]
			(M) { \lhm{$\Sigma_\theta$}}; &
    \node		(node)  {}; \\
    \node[draw, circle, minimum size=0.05cm]
			(sum_2) {};
			\node at (sum_2.north east)[right=-3pt,yshift=+3pt]{ $+$};
                \node at (sum_2.south east)[right=-3pt,yshift=-3pt]{$+$}; &
    \node[draw, rectangle, minimum width=0.75cm, minimum height=0.5cm]
			(F) { $\nabla F(\cdot)$}; &
    \node	(sum_3) {}; \\
    \node[draw, rectangle, minimum width=0.75cm, minimum height=0.5cm, xshift=0.75cm]
			(nablagx) { \lhg{$\nabla G(\theta) $}}; &
    \node[draw, rectangle, minimum width=0.75cm, minimum height=0.75cm]
			(G) {\lhm{$\Sigma_\lambda^+$}}; &
    \node[draw, rectangle, minimum width=0.75cm,
    minimum height=0.5cm,xshift=-0.75cm]
			(gx) { \lhg{$G(\cdot)$}}; \\
};

\node[draw, rectangle,
			below=of G,
			yshift=+0.5cm,
			minimum width=0.75cm, 	minimum height=0.75cm]
			(H) { \lhm{$\Sigma_\mu$}};

\node[draw, circle,
			minimum size=0.6cm]
			(B) at (0,-0.75) {};
			\node at (B.center){ $\mathcal{A}$};
\node[draw, circle,
			right=of B,
			xshift= +5.8cm,
			minimum size=0.6cm]
			(Bt)   {};
			\node at (Bt.center){$\mathcal{A}^T$};
 
\node[draw,thick,dotted,
minimum width =5.2cm, minimum height =3.1cm]
(box) at (3.8,0) {};
\node at (box.north east)[left=8pt,yshift=7pt,text width=3cm]{Node dynamics};
 \node[draw,thick,dotted,
minimum width =5.2cm, minimum height =1.2cm]
(box) at (3.8,-2.35) {};
\node[left=8pt,yshift=-7pt,text width=3cm] at (box.south east){Edge dynamics};

\draw[-stealth] (sum_1.east) -- (M.west)
   node[midway,above]{ $\phi$};
\draw (M.east) -- (node.center)
	node[pos=0.5,above]{ $\theta$};
\draw(node.center) -- (sum_3.center);
\draw[-stealth] (sum_2.north) -- (sum_1.south);
\draw[-stealth] (sum_3.center) -- (F.east);
\draw[-stealth] (F.west) -- (sum_2.east);
\draw[-stealth] (sum_3.center) |- (gx.east);
\draw[-stealth] (nablagx.west) -| (sum_2.south)
    node[pos=0.75,left]{ $\eta$};
\draw[-stealth] (gx.west) -- (G.east);
\draw[-stealth] (G.west) -- (nablagx.east)
    node[midway,above]{$\lambda$};
\draw[-stealth] (node.center) -| (Bt.north);
\draw[-stealth] (B.north) |- (sum_1.west)
    node[pos=0.75,above]{\lhf{$\psi$}};
\draw[-stealth] (Bt.south) |- (H.east)
    node[pos=0.75,above]{ $\omega$};
\draw[-stealth] (H.west) -| (B.south)
    node[pos=0.25,above]{ $\mu$};

\end{tikzpicture}

%% file: Images/tikz_network.tex
\begin{tikzpicture}
\node[draw, circle,
			minimum size=0.1cm,
			align=center]
			(nu1) at (0,0) {$\nu_1$};
\node[draw, circle,
			minimum size=0.1cm,
			align=center]
			(nu2) at (-1,-1) {$\nu_2$};
\node[draw, circle,
			minimum size=0.1cm,
			align=center]
			(nu3) at (1,-1) {$\nu_3$};

\draw[-stealth] (nu1.south west) -- (nu2.north east);
\draw[-stealth] (nu1.south east) -- (nu3.north west);
\draw[-stealth] (nu2.east) -- (nu3.west);

\end{tikzpicture}

%% file: Images/tikz_inequality.tex
\begin{tikzpicture}


\node[matrix,column sep=1cm, row sep =0.25cm] at (3.75,0)
{
    \node[draw, circle, minimum size=0.05cm](sum_1){};		
			\node at (sum_1.north west)[left=-2pt,yshift=+3pt]{ $-$};
			\node at (sum_1.south east)[right=-2pt,yshift=-3pt]{$-$}; &
    \node[draw, rectangle,
			minimum width=0.75cm, 	minimum height=0.75cm]
			(M) { \lhm{$\Sigma_\theta$}}; &
    \node		(node)  {}; \\
    \node (sum_2) {}; &
    \node[draw, rectangle, minimum width=0.75cm, minimum height=0.5cm]
			(F) { $\nabla F(\cdot)$}; &
    \node	(sum_3) {}; \\
};

\node[draw, rectangle,
			below=of F,
			yshift=+0.5cm,
			minimum width=0.75cm, 	minimum height=0.75cm]
			(G) { \lhm{$\Sigma_\lambda^+$}};
\node[draw, rectangle, minimum width=0.75cm, minimum height=0.5cm, left=of G] (nablagx) { \lhg{$\nabla H(\omega) $}};
\node[draw, rectangle, minimum width=0.75cm,
    minimum height=0.5cm,right=of G]
			(gx) { \lhg{$H(\cdot)$}};

\node[draw, circle,
			minimum size=0.6cm]
			(B) at (0,-0.75) {};
			\node at (B.center){ $\mathcal{R}^T$};
\node[draw, circle,
			right=of B,
			xshift= +5.8cm,
			minimum size=0.6cm]
			(Bt)   {};
			\node at (Bt.center){$\mathcal{R}$};
 
\node[draw,thick,dotted,
minimum width =5.1cm, minimum height =2cm]
(box) at (3.9,0) {};
\node at (box.north east)[right=3pt,yshift=-3pt]{Nodes};
 \node[draw,thick,dotted,
minimum width =5.1cm, minimum height =1cm]
(box) at (3.9,-1.7) {};
\node at (box.south east)[right=3pt,yshift=6pt]{Edges};

\draw[-stealth] (sum_1.east) -- (M.west)
   node[midway,above]{ $\phi$};
\draw (M.east) -- (node.center)
	node[pos=0.5,above]{ $\theta$};
\draw(node.center) -- (sum_3.center);
\draw[-stealth] (sum_2.center) -- (sum_1.south);
\draw[-stealth] (sum_3.center) -- (F.east);
\draw (F.west) -- (sum_2.center);
\draw[-stealth] (node.center) -| (Bt.north);
\draw[-stealth] (B.north) |- (sum_1.west)
    node[pos=0.75,above]{\lhf{$\eta$}};
\draw[-stealth] (Bt.south) |- (gx.east)
    node[pos=0.75,above]{ $\omega$};
\draw[-stealth] (nablagx.west) -| (B.south);
\draw[-stealth] (gx.west) -- (G.east);
\draw[-stealth] (G.west) -- (nablagx.east)
    node[midway,above]{$\lambda$};

\end{tikzpicture}

%% file: Images/tikz_feed_forward.tex
\begin{tikzpicture}
\small

\node[draw, circle,
            minimum size=0.5cm]
            (neg) at (0,0) {};
            \node at (neg.center){  $-1$};
\node[draw, circle,
            minimum size=0.05cm,
            right=of neg,
            xshift=0cm]
            (sum_1) {};
            \node at (sum_1.north)[above=-2pt]{  $+$};
            \node at (sum_1.south east)[right=-2pt]{ $-$};
\node[draw, rectangle,
            right=of sum_1,
                   xshift=-0.25cm,
            minimum width=0.5cm, 	minimum height=0.5cm]
            (M) {  \lhm{$\Sigma_\theta$}};
\node[right=of M, xshift=-0.25cm]
            (node_top)  {};
\node[draw, rectangle,
            below=of M,
            yshift = +0.75cm,
            minimum width=0.75cm, 	minimum height=0.5cm]
            (F) { $\nabla F(\cdot) $};
\node[draw, rectangle,
            below=of F,
            yshift=+0.5cm,
            minimum width=0.5cm, 	minimum height=0.5cm]
            (H) {  \lhm{$\tilde{\Sigma}_\mu$}};
\node[draw, rectangle,
            below=of H,
            yshift=+0.75cm,
            minimum width=0.3cm, 	minimum height=0.3cm]
            (d) { $\check{d}$};
\node[right=of H, xshift=-0.25cm]
            (node_bottom)  {};
\node[draw, circle,
            minimum size=0.05cm,
            left=of H,
                xshift=0.25cm]
            (sum_2) {};
            \node at (sum_2.east)[above=2pt]{  $+$};
            \node at (sum_2.south)[left=2pt]{  $+$};
\node[draw, circle,
            minimum size=0.5cm]
            (B) at (0,-1) {};
            \node at (B.center){  $\mathcal{A}$};
\node[draw, circle,
            right=of M,
            xshift= 0.4cm,
            yshift= -1cm,
            minimum size=0.5cm]
            (Bt)   {};
            \node at (Bt.center){ $\mathcal{A}^T$};
            
\draw[-stealth] (neg.east) -- (sum_1.west)
    node[midway,above]{ $-\psi$};
\draw[-stealth] (F.west) -| (sum_1.south);
\draw[-stealth] (sum_1.east) -- (M.west)
   node[midway,above]{ $\phi$};
\draw (M.east) -- (node_top.center);
\draw[-stealth] (node_top.center) |- (F.east);
\draw[-stealth] (node_top.center) -| (Bt.north)
    node[pos=0.25,above]{ $\theta$};
\draw (Bt.south) |- (node_bottom.center)
    node[pos=0.75,above]{ $\omega$};
\draw[-stealth] (node_bottom.center) -- (H.east);
\draw[-stealth] (node_bottom.center) |- (d.east);
\draw[-stealth] (H.west) -- (sum_2.east)
    node[midway,above]{$\tilde{\mu}$};
\draw[-stealth] (d.west) -| (sum_2.south);
\draw[-stealth] (sum_2.west) -| (B.south)
    node[pos=0.25,above]{ $\mu$};
\draw[-stealth] (B.north) -- (neg.south);

\end{tikzpicture}

%% file: Images/tikz_feedback.tex
\begin{tikzpicture}
\small

\node[draw, circle,
            minimum size=0.5cm]
            (neg) at (0,0) {};
            \node at (neg.center){ $-1$};
\node[draw, circle,
            minimum size=0.05cm,
            right=of neg,
            xshift=-0.2cm]
            (sum_1) {};
            \node at (sum_1.north)[above=-2pt]{ $+$};
            \node at (sum_1.south east)[right=-2pt,yshift=-3pt]{$-$};
\node[draw, circle,
            minimum size=0.05cm,
            right=of sum_1,
            xshift=-0.4cm]
            (sum_2) {};
            \node at (sum_2.north east)[above=-2pt,xshift=3pt]{ $-$};
            \node at (sum_2.south west)[left=-3pt,yshift=-3pt]{$+$};
\node[draw, rectangle,
            right=of sum_2,
            xshift=-0.4cm,
            minimum width=0.5cm, 	minimum height=0.5cm]
            (M) { \lhm{$\Sigma_\theta$}};
\node[right=of M,
            xshift=-0.7cm]
            (node_right_1)  {};
\node[right=of node_right_1,
            xshift=-0.7cm]
            (node_right_2)  {};
\node[draw, circle,
                above=of sum_2,
                yshift=-0.75cm,
            minimum size=0.5cm]
            (B_top) {};
            \node at (B_top.center){ $\mathcal{A}$};
\node[draw, circle,
                above=of node_right_1,
                yshift=-0.75cm,
            minimum size=0.5cm]
            (Bt_top) {};
            \node at (Bt_top.center){ $\mathcal{A}^T$};
\node[draw, rectangle,
            below=of M,
            yshift = +0.75cm,
            minimum width=0.5cm, 	minimum height=0.5cm]
            (F) { $\nabla F(\cdot) $};
\node[draw, rectangle,
            above=of M,
             yshift=-0.25cm,
            minimum width=0.4cm, 	minimum height=0.3cm]
            (d) {$\check{d}$};
\node[draw, rectangle,
            below=of F,
            yshift=+0.5cm,
            minimum width=0.5cm, 	minimum height=0.5cm]
            (H) { \lhm{$\tilde{\Sigma}_\mu$}};
\node[draw, circle,
            minimum size=0.5cm]
            (B) at (0,-1) {};
            \node at (B.center){ $\mathcal{A}$};
\node[draw, circle,
            right=of M,
            xshift= 0.4cm,
            yshift= -1cm,
            minimum size=0.5cm]
            (Bt)   {};
            \node at (Bt.center){$\mathcal{A}^T$};
            
\draw[-stealth] (neg.east) -- (sum_1.west)
    node[midway,above]{ $-\tilde{\psi}$};
\draw[-stealth] (F.west) -| (sum_1.south);
\draw[-stealth] (sum_1.east) -- (sum_2.west)
   node[midway,above]{ $\tilde{\phi}$};
\draw[-stealth] (sum_2.east) -- (M.west)
    node[midway,above]{ $\phi$};
\draw (M.east) -- (node_right_1.center);
\draw[-stealth] (node_right_1.center) -- (Bt_top.south);
\draw[-stealth] (Bt_top.north) |- (d.east);
\draw[-stealth] (d.west) -| (B_top.north);
\draw[-stealth] (B_top.south) -- (sum_2.north);
\draw[-stealth] (node_right_2.center) |- (F.east);
\draw[-stealth] (node_right_1.center) -| (Bt.north)
    node[pos=0.25,above]{ $\theta$};
\draw (Bt.south) |- (H.east)
    node[pos=0.75,above]{ $\omega$};
\draw[-stealth] (H.west) -| (B.south)
    node[pos=0.25,above]{ $\tilde{\mu}$};
\draw[-stealth] (B.north) -- (neg.south);

\end{tikzpicture}